\documentclass{gtart_h}  

\def\ifplaintex{\expandafter\ifx\csname documentclass\endcsname\relax}

\ifplaintex 
\else
\fi

\expandafter\ifx\csname beginpicture\endcsname\relax
\expandafter\ifx\csname documentclass\endcsname\relax
\input pictex \else
\input prepictex \input pictex \input postpictex \fi\fi

\def\gt{{\mathsurround=0pt\it $\cal G\mskip-2mu$eometry \&\ 
$\cal T\!\!$opology}}        

\def\gtp{{\mathsurround=0pt\it $\cal G\mskip-2mu$eometry \&\ 
$\cal T\!\!$opology $\cal P\!$ublications}}  

\def\lognumber#1{\def\thelognumber{#1}}
\def\volumenumber#1{\def\thevolumenumber{#1}}
\def\papernumber#1{\def\thepapernumber{#1}}
\def\volumeyear#1{\def\thevolumeyear{#1}}

\def\pagenumbers#1#2{\def\startpage{#1}\def\finishpage{#2}}
\def\published#1{\def\publishdate{#1}}
\def\proposed#1{\def\theproposer{#1}}
\def\seconded#1{\def\theseconders{#1}}
\def\received#1{\def\receiveddate{#1}}
\def\revised#1{\def\reviseddate{#1}}
\def\accepted#1{\def\accepteddate{#1}}
\def\asciititle#1{\def\theasciititle{#1}}
\def\covertitle#1{\def\thecovertitle{#1}}

\def\asciiaddress#1{\def\theasciiaddress{#1}}

\long\def\asciiabstract#1{\long\def\theasciiabstract{#1}}
\def\asciikeywords#1{\def\theasciikeywords{#1}}

\let\\\par\let\thelognumber\relax
\let\thevolumenumber\relax\let\thepapernumber\relax
\let\thevolumeyear\relax\let\thesamplenumber\relax\let\startpage\relax
\let\finishpage\relax\let\publishdate\relax\let\receiveddate\relax
\let\reviseddate\relax\let\accepteddate\relax\let\theasciititle\relax
\let\thecovertitle\relax\let\theasciiauthors\relax\let\theasciiaddress\relax
\let\theasciiabstract\relax\let\theasciikeywords\relax
\let\theasciiemail\relax\let\theshortauthors\relax\let\theshorttitle\relax

\long\def\maketitlep{   

\count0=\startpage

\gt\hfill      
\beginpicture
\setcoordinatesystem units <0.33truein, 0.33truein> point at 2.2 0.9
\setplotsymbol ({$\cal G$})
\plotsymbolspacing=9truept
\circulararc 315 degrees from 0 1 center at 0 0
\setplotsymbol ({$\cal T$})
\circulararc 315 degrees from 1 -1 center at 1 0
\endpicture
\break
{\small\ifx\thesamplenumber\relax 
Volume \else Sample
\fi\thevolumenumber\ (\thevolumeyear)
\startpage--\finishpage\nl
Published: \publishdate}
\vglue 0.5truein plus 0.4fil minus 0.1truein

{\parskip=0pt\leftskip 0pt plus 1fil\def\\{\par\smallskip}{\ifplaintex\large
\else\Large\fi\bf\thetitle}\par\medskip}   

\vglue 0pt plus 0.1fil 

{\parskip=0pt\leftskip 0pt plus 1fil\def\\{\par}{\sc\theauthors}
\par\medskip}

\vglue 0pt plus 0.1fil 

{\small\parskip=0pt\let\newline\\
{\leftskip 0pt plus 1fil\def\\{\par}{\sl\theaddress}\par}
\expandafter\ifx\theemail\relax    
\relax\else\vglue 5pt plus 0.02fil minus 2pt\def\\{\stdspace{\rm 
and}\stdspace} 
\cl{Email:\stdspace\tt\theemail}\fi
\ifx\theurl\relax                  
\relax\else\vglue 5pt plus 0.02fil minus 2pt\def\\{\stdspace{\rm 
and}\stdspace}
\cl{URL:\stdspace\tt\theurl}\fi\par}

\vglue 7pt plus 0.3fil minus 3pt

{\bf Abstract}
\vglue 5pt plus 0.1fil minus 2pt

\theabstract

\vglue 7pt plus 0.3fil minus 3pt

{\bf AMS Classification numbers}\quad Primary:\quad \theprimaryclass

Secondary:\quad \thesecondaryclass

\vglue 5pt plus 0.3fil minus 2pt

{\bf Keywords:}\quad \thekeywords

\vglue 10pt plus 0.5fil minus 5pt

{\small  Proposed: \theproposer\hfill Received: \receiveddate\nl
Seconded: \theseconders\hfill 
\ifx\reviseddate\relax                         
Accepted: \accepteddate                        
\else
Revised: \reviseddate                          
\fi}
\eject
}       

\let\maketitlepage\maketitlep
\let\maketitle\maketitlepage

\font\phead=cmsl9 scaled 950
\font=cmsl9 scaled 1050
\font\pnum=cmbx10 scaled 913
\font\lnum=cmbx10 
\font\pfoot=cmsl9 scaled 950
\font=cmsl9 scaled 1050
\ifplaintex
\headline{\vbox to 0pt{\vskip -4.5mm\line{\small\phead\ifnum
\count0=\startpage ISSN 1364-0380 (on line)
1465-3060 (printed) \hfill {\pnum\folio}\else\ifodd\count0\def\\{ }%
\ifx\theshorttitle\relax\thetitle\else\theshorttitle\fi\hfill{\pnum\folio}
\else\def\\{ and }{\pnum\folio}\hfill\ifx\theshortauthors\relax\theauthors
\else\theshortauthors\fi\fi\fi}\vss}}
\footline{\vbox to 0pt{\vglue 0mm\line{\small\pfoot\ifnum\count0=\startpage
\copyright\ \gtp\hfill\else
\gt, Volume \thevolumenumber\ (\thevolumeyear)\hfill\fi}\vss
}}
\else
\makeatletter
\def\@oddhead{{\small\lhead\ifnum\count0=\startpage ISSN 1364-0380 (on line)
1465-3060 (printed) \hfill {\lnum\number\count0}\else\ifodd\count0
\def\\{ }\ifx\theshorttitle\relax \thetitle \else\theshorttitle\fi\hfill
{\lnum\number\count0}\else\def\\{ and }{\lnum\number\count0}
\hfill\ifx\theshortauthors\relax 
\theauthors\else\theshortauthors\fi\fi\fi}}\def\@evenhead{\@oddhead}
\def\@oddfoot{\small\lfoot\ifnum\count0=\startpage\copyright\ \gtp\hfill\else
\gt, Volume \thevolumenumber\ (\thevolumeyear)\hfill\fi}
\def\@evenfoot{\@oddfoot}
\makeatother
\fi

\newwrite\gtoutfile
\long\gdef\makeheadfile{  
{\def\\{, }\def\s{ }
\immediate\openout\gtoutfile head.xxx
\immediate\write\gtoutfile{Proxy-for: \ifx\theasciiauthors\relax
\theauthors\else\theasciiauthors\fi\s<\ifx\theasciiemail\relax\theemail\else\theasciiemail\fi>}
\immediate\write\gtoutfile{\noexpand\\}
\immediate\write\gtoutfile{Authors: \ifx\theasciiauthors\relax
\theauthors\else\theasciiauthors\fi}
{\def\\{ }\immediate\write\gtoutfile{Title: \ifx\theasciititle\relax
\thetitle\else\theasciititle\fi}}
\immediate\write\gtoutfile{Subj-class: GT or SG or MG etc}
\immediate\write\gtoutfile{MSC-class: \theprimaryclass\ifx\thesecondaryclass\relax\else, \thesecondaryclass\fi}
\immediate\write\gtoutfile{Journal-ref: Geom. Topol. \thevolumenumber
(\thevolumeyear) \startpage-\finishpage}
\immediate\write\gtoutfile{Comments: Published by Geometry and Topology at}
\immediate\write\gtoutfile{\s\s http://www.maths.warwick.ac.uk/gt/GTVol\thevolumenumber/paper\thepapernumber.abs.html}
\immediate\write\gtoutfile{\noexpand\\}
\immediate\write\gtoutfile{}
\ifx\theasciiabstract\relax
\immediate\write\gtoutfile{\theabstract}\else
\immediate\write\gtoutfile{\theasciiabstract}\fi
\immediate\write\gtoutfile{}
\immediate\write\gtoutfile{\noexpand\\}
\immediate\write\gtoutfile{}
\immediate\closeout\gtoutfile}}  

\def\maketitlepage{\maketitlep\makeheadfile}
\let\maketitle\maketitlepage

\lognumber{370}
\received{14 October 2003}
\volumenumber{8}\papernumber{39}\volumeyear{2004}
\pagenumbers{1427}{1470}
\revised{26 November 2004}
\published{27 November 2004}
\accepted{29 September 2004}
\proposed{Martin Bridson}
\seconded{Benson Farb, Walter Neumann}

\usepackage{amsmath,amssymb,rotating,mdwlist,export}

\newtheorem{thm}{Theorem}[section]

\newtheorem*{thm*}{Theorem}
\newtheorem{dfn}[thm]{Definition} 
\newtheorem*{dfn*}{Definition}

\newtheorem{cor}[thm]{Corollary}
\newtheorem*{cor*}{Corollary}

\newtheorem{prop}[thm]{Proposition} 
\newtheorem*{prop*}{Proposition} 
 
\newtheorem{lem}[thm]{Lemma} 
\newtheorem*{lem*}{Lemma}

\newtheorem*{claim*}{Claim} 
 
\newtheorem*{fact*}{Fact} 
\newtheorem{fact}[thm]{Fact}

\theoremstyle{remark}
\newtheorem*{rem*}{Remark}
\newtheorem*{example*}{Example}

\newenvironment{SauveCompteurs}[1]{%
\newcommand{\monparametre}{#1}
\openexport{\monparametre_sauve}
  \Export{thm}\Export{section}\Export{subsection}\Export{subsubsection}
\closeexport}{}

\newenvironment{UtiliseCompteurs}[1]{%
\newcommand{\monparametre}{#1}
\openexport{\monparametre_aux}
  \Export{thm}\Export{section}\Export{subsection}\Export{subsubsection}
\closeexport
\PackageInfo{export}{\MessageBreak
Importations from \monparametre_sauve.xpt\MessageBreak}%
\InputIfFileExists{\monparametre_sauve.xpt}{\relax}{\relax}%
\renewcommand{\label}[1]{}
}{%
\PackageInfo{export}{\MessageBreak
Importations from \monparametre_aux.xpt\MessageBreak}%
\InputIfFileExists{\monparametre_aux.xpt}{\relax}{\relax}}

\newlength{\espaceavantspecialthm}
\newlength{\espaceapresspecialthm}
{\normalfont \vskip \espaceapresspecialthm}

\newenvironment{specialthm*}[1]{
\vskip\espaceavantspecialthm \noindent \textbf{#1} \itshape}%
{\normalfont \vskip \espaceapresspecialthm}

\newlength{\espaceavantenonce}
\newlength{\espaceapresenonce}
\newcommand{\fontetitreun}[1]{\textbf{#1}} 
\newcommand{\fontetitredeux}[1]{\textit{#1}} 

{\normalfont \vskip \espaceapresenonce}

\newenvironment{enonce1*}[1]{
\vskip\espaceavantenonce \noindent \fontetitreun{#1} \itshape}%
{\normalfont \vskip \espaceapresenonce}

{\vskip \espaceapresenonce}

\newenvironment{enonce2*}[1]{
\vskip\espaceavantenonce \noindent \fontetitredeux{#1} }%
{\vskip \espaceapresenonce}

\newcommand{\es}{\emptyset}
\renewcommand{\phi}{\varphi} 
\newcommand{\m} {^{-1}} 
\newcommand{\eps} {\varepsilon}

\newcommand {\ra} {\rightarrow}

\newcommand {\onto} {\twoheadrightarrow}

\newcommand{\actson}{\,\raisebox{1.8ex}[0pt][0pt]{\begin{turn}{-90}\ensuremath{\circlearrowright}\end{turn}}\,}

\newcommand{\ol}[1]{\overline{#1}}

\newcommand{\tensor}{\otimes}
\newcommand{\dunion}{\sqcup}

\newcommand{\ie} {i.e.\ } 
\newcommand{\Tmin}{T_{\mathrm{min}}}

\newcommand {\calb} {{\mathcal {B}}}

\newcommand {\calf} {{\mathcal {F}}}   
\newcommand {\calg} {{\mathcal {G}}}

\newcommand {\calt} {{\mathcal {T}}}

\newcommand{\bbZ} {{\mathbb{Z}}}   
\newcommand{\bbN} {{\mathbb{N}}}   
\newcommand{\bbR} {{\mathbb{R}}}   
\newcommand{\bbQ} {{\mathbb{Q}}}

\DeclareMathOperator{\Tr}{Tr}
\DeclareMathOperator{\Isom}{Isom}
\DeclareMathOperator{\Rk}{Rk}
\DeclareMathOperator{\Stab}{Stab}  
\DeclareMathOperator*{\Fix}{Fix}

\newcommand{\Id} {\mathrm{Id}}

\newcommand{\eqv} {{\sim}}

\let \para \S
\renewcommand{\S} {\Sigma}

\newcommand{\olY}{{\ol{Y}}\!}

\newcommand{\olZ}{{\ol{Z}}}
\newcommand{\olT}{{\ol{T}}{}}

\newcommand{\olG}{{\ol{\Gamma}}}
\newcommand{\calbt}{{\calb\calt}}

\begin{document}

\title{Limit groups and groups acting freely on $\bbR^n$--trees}
\covertitle{Limit groups and groups acting freely on ${\noexpand\bf R}^n$--trees}
\asciititle{Limit groups and groups acting freely on R^n-trees}
\author{Vincent Guirardel}

\address{Laboratoire E. Picard, UMR 5580, B\^at 1R2\\
Universit\'e Paul Sabatier, 118 rte de Narbonne\\
31062 Toulouse cedex 4, France}
\email{guirardel@picard.ups-tlse.fr}
\asciiaddress{Laboratoire E. Picard, UMR 5580, Bat 1R2\\
Universite Paul Sabatier, 118 rte de Narbonne\\
31062 Toulouse cedex 4, France}

\asciikeywords{R^n-tree, limit group, finite presentation}

\begin{abstract}
We give a simple proof of the finite presentation of Sela's limit groups
by using free actions on $\bbR^n$--trees.
We first prove that Sela's limit groups do have a free action on an $\bbR^n$--tree.
We then prove that a finitely generated group having a free action on an $\bbR^n$--tree can be
 obtained from free abelian groups and surface groups by a finite sequence 
of free products and amalgamations over cyclic groups.
As a corollary, such a group is finitely presented, has a finite classifying space,
its abelian subgroups are finitely generated and contains only
finitely many conjugacy classes of non-cyclic maximal abelian subgroups.
\end{abstract}

\asciiabstract{%
We give a simple proof of the finite presentation of Sela's limit
groups by using free actions on R^n-trees.  We first prove that
Sela's limit groups do have a free action on an R^n-tree.  We
then prove that a finitely generated group having a free action on an
R^n-tree can be obtained from free abelian groups and surface
groups by a finite sequence of free products and amalgamations over
cyclic groups.  As a corollary, such a group is finitely presented,
has a finite classifying space, its abelian subgroups are finitely
generated and contains only finitely many conjugacy classes of
non-cyclic maximal abelian subgroups.}

\maketitle

\section{Introduction}
Limit groups have been introduced by Z Sela in the first paper of his
solution of Tarski's problem \cite{Sela_diophantine1}. These groups appeared
to coincide with the long-studied class of finitely generated fully
residually free groups
(see \cite{Baumslag_residually}, \cite{Baumslag_generalised},
\cite{KhMy_irreducible1,KhMy_irreducible2}, 
\cite{Chi_introduction} and references).

A limit group is a limit of free groups in the space of marked groups.
More precisely, if $n$ is a fixed integer, a marked group is a group
together with an ordered generating family $S=(s_1,\dots,s_n)$.
Two marked groups $(\Gamma,S)$ and $(\Gamma',S')$ are close to each other in this topology if 
for some large $R$, $(\Gamma,S)$ and $(\Gamma',S')$ have exactly the same relations of length at most $R$
(see section \ref{sec_prelim_limit}).

Limit groups  have several equivalent characterisations: a finitely generated group $G$ is a limit group
if and only if it is fully residually free,
 if and only if it has the same universal theory as a free group,
if and only if it is a subgroup of a non-standard free group \cite{Remeslennikov_exist_siberian,CG_compactifying}.
 
One of the main results about limit groups is a structure theorem
due to Kharlampovich--Myasnikov, Pfander and Sela \cite{KhMy_irreducible1,KhMy_irreducible2,Pfander_finitely,Sela_diophantine1}.
This theorem claims that a limit group can be inductively obtained from free abelian groups
and surface groups by taking free products and amalgamations over $\bbZ$ (see Theorem \ref{devissage_simple} below).
This structure theorem implies that a limit group is finitely presented, and that its abelian subgroups are finitely generated.
The goal of the paper is to give a simpler  proof of this result in the broader context of groups acting freely 
on $\bbR^n$--trees.

After completing this work, the author learnt about the unpublished thesis of Shalom Gross, a student of
Z~Sela,  proving the finite presentation of finitely generated groups having a free action on an $\bbR^n$--tree \cite{Gross_these}. 
Both proofs deeply rely on Sela's structure theorem for super-stable actions of finitely generated groups
on $\bbR$--trees (\cite[Theorem 3.1]{Sela_acylindrical}, see also Theorem \ref{thm_structure} below).
However, Gross does not state a d\'evissage theorem over cyclic groups,
but over finitely generated abelian groups.

Let's recall briefly the definition of a $\Lambda$--tree. Given a totally ordered abelian group $\Lambda$,
there is a natural notion of $\Lambda$--metric space where the distance function takes its values in $\Lambda$.
If $\Lambda$ is archimedean, then $\Lambda$ is isomorphic to a subgroup of $\bbR$ and we have a metric in the usual sense.
When $\Lambda$ is not archimedean, there are elements which are infinitely small compared to other elements.
A typical example is when $\Lambda=\bbR^n$ endowed with the lexicographic ordering.

A $\Lambda$--tree may be defined as a geodesic $0$--hyperbolic $\Lambda$--metric space.
Roughly speaking, an $\bbR^n$--tree may be thought of as a kind of bundle over an $\bbR$--tree
where the fibres are (infinitesimal) $\bbR^{n-1}$--trees.

In his list of research problems, Sela conjectures that a finitely generated group is a limit group
if and only if it acts freely on an $\bbR^n$--tree \cite{Sela_problems}.
However, it is known that the fundamental group $\Gamma=\langle a,b,c\,|\, a^2b^2c^2=1\rangle$
of the non-orientable surface $\Sigma$ of Euler characteristic $-1$
is not a limit group since three elements in a free group satisfying $a^2b^2c^2=1$ must commute 
(\cite{Lyndon_equation},\cite[p.249]{Chi_book}).
But this group acts freely on a $\bbZ^2$--tree: $\Sigma$ can be obtained
by gluing together the two boundary components of a twice punctured projective plane,
so $\Gamma$ can be written as an HNN extension $F_2*_\bbZ$. The $\bbZ^2$--tree can be roughly described
as the Bass--Serre tree of this HNN extension, but where one blows up each vertex into
an infinitesimal tree corresponding to a Cayley graph of $F_2$ (see \cite[p.237]{Chi_book} for details).

In this paper, we start by giving a proof that every limit group acts freely on an $\bbR^n$--tree. 
This is an adaptation a theorem by Remeslennikov saying that a fully residually free
groups act freely on a $\Lambda$--tree where $\Lambda$ has finite $\bbQ$--rank, \ie $\Lambda\tensor\bbQ$ is finite dimensional
(\cite{Remeslennikov_exist_ukrainian}, see also \cite[Theorem  5.10]{Chi_book}).
However, Remeslennikov claims that $\Lambda$ can be chosen finitely generated, but this relies
on a misquoted result about valuations (see section \ref{sec_remeslennikov}).

We actually prove that there is a closed subspace of the space of marked groups
consisting of groups acting freely on $\bbR^n$--trees. In the following statement, 
an action of a group $\Gamma$ on a Bruhat--Tits tree
is the action on the Bruhat--Tits tree of $SL_2(K)$ induced by a morphism $j\co \Gamma\ra SL_2(K)$ 
where $K$ is a valuation field. Note that $K$  may vary with $\Gamma$.

\newcommand{\thmclosed}{%
\begin{thm}[Acting freely on Bruhat--Tits trees is closed, see also \cite{Remeslennikov_exist_ukrainian}]
Let $\calbt\subset\calg_n$ be the set of marked groups $(\Gamma,S)$ 
having a free action on a Bruhat--Tits tree.

Then $\calbt$ is closed in $\calg_n$, and $\calbt$ consists in groups acting freely on $\Lambda$--trees
where $\Lambda\tensor\bbQ$ has dimension at most $3n+1$ over $\bbQ$. 
In particular, $\calbt$ consists in groups acting freely on $\bbR^{3n+1}$--trees where $\bbR^{3n+1}$ has the
lexicographic ordering.
\end{thm}}
\begin{UtiliseCompteurs}{thm_closed}\renewcommand{\label}[1]{}
\thmclosed
\end{UtiliseCompteurs}

\newcommand{\corRemeslennikov}{%
\begin{cor}{\rm\cite{Remeslennikov_exist_ukrainian}}\qua
A limit group has a free action on an $\bbR^n$--tree.
\end{cor}}
\begin{UtiliseCompteurs}{Remeslennikov}\corRemeslennikov
\end{UtiliseCompteurs}

\begin{rem*}
As a corollary of their study of the structure of limit groups \cite{KhMy_irreducible1,KhMy_irreducible2},
  Kharlampovich and Myasnikov prove the more precise result that a limit group 
is a subgroup of an iterated free extension of centralizers of a free group, 
and has therefore a free action on a $\bbZ^n$--tree \cite[Corollary 6]{KhMy_irreducible2}.
An alternative proof of this fact using Sela's techniques is given in \cite{CG_compactifying}.
\end{rem*}

The main result of the paper is the following structure theorem for groups acting freely on $\bbR^n$--trees
(see theorem \ref{cyclic_devissage} for a more detailed version).
In view of the previous corollary, this theorem applies to limit groups.

\newcommand{\devissageSimple}{%
\begin{thm}[D\'evissage theorem, simple version. See also {\cite[Corollary~6.6]{Gross_these}}]\label{devissage_simple}
Consider a finitely generated, freely indecomposable group $\Gamma$ having a free action on an $\bbR^n$--tree.
Then $\Gamma$ can be written as the fundamental group of a finite graph of groups 
with cyclic edge groups and where each vertex group is finitely generated and has
a free action on an $\bbR^{n-1}$--tree.
\end{thm}}
\begin{UtiliseCompteurs}{devissage_simple}
\devissageSimple
\end{UtiliseCompteurs}

For $n=1$, Rips theorem says that $\Gamma$ (which is supposed to be freely indecomposable)
is either a free abelian group, or a surface group (see \cite{GLP1,BF_stable}).
Hence, a limit group can be obtained from abelian and surface groups by a finite sequence
of free products and amalgamations over $\bbZ$. It is therefore easy to deduce the following result:

\newcommand{\corFP}{%
\begin{cor}[See also \cite{Gross_these}]
  Let $\Gamma$ be a finitely generated group having a free action on an $\bbR^n$--tree.
Then 
\begin{itemize*}
\item $\Gamma$ is finitely presented \cite[Corollary 6.6]{Gross_these};
\item if $\Gamma$ is not cyclic, then its first Betti number is at least $2$;
\item there are finitely many conjugacy classes of non-cyclic maximal abelian subgroups in $\Gamma$,
and abelian subgroups of $\Gamma$ are finitely generated. More precisely,
one has the following bound on the ranks of maximal abelian subgroups:
$$\sum_A (\Rk A-1)\leq b_1(\Gamma)-1$$
where the sum is taken over the set of conjugacy classes of non-cyclic maximal abelian subgroups of $\Gamma$,
and where $b_1(\Gamma)$ denotes the first Betti number of $\Gamma$;
\item $\Gamma$ has a finite classifying space, and the cohomological dimension of $\Gamma$ is at most $\max(2,r)$ where $r$ is the maximal rank
of an abelian subgroup of $\Gamma$.
\end{itemize*}
\end{cor}}
\begin{UtiliseCompteurs}{cor_FP}%
\corFP
\end{UtiliseCompteurs}

\begin{rem*}
  A combination theorem by Dahmani also shows that $\Gamma$ is hyperbolic relative to its non-cyclic abelian subgroups 
\cite{Dahmani_combination}.
\end{rem*}

\begin{cor*}{\rm\cite{Sela_diophantine1,KhMy_irreducible1,KhMy_irreducible2,Pfander_finitely}}\qua
A limit group is finitely presented, its abelian subgroups are finitely generated, 
it has only finitely many conjugacy classes of maximal non-cyclic abelian subgroups, and it has
a finite classifying space.
\end{cor*}

Finally, we can also easily derive from the d\'evissage theorem 
the existence of a \emph{principal} splitting, a major step in Sela's proof
of the finite presentation of limit groups (see corollary \ref{cor_principal} and \cite[Theorem 3.2]{Sela_diophantine1}). 

Unlike Sela's proof, the proof we give doesn't need any JSJ theory,
and does not use the shortening argument.
The proof is also much shorter than the one by Kharlampovich--Myasnikov
in \cite{KhMy_irreducible1,KhMy_irreducible2} using algebraic geometry over groups, and
the study of equations in free groups. 

The paper is organized as follows: after some premilinaries in section \ref{sec_prelim},
section \ref{sec_remeslennikov} is devoted to the proof of the fact that limit groups act freely on $\bbR^n$--trees.
Section \ref{sec_gluing} sets up some preliminary work on graph of actions on $\Lambda$--trees, which encode
how to glue equivariantly some $\Lambda$--trees to get a new $\Lambda$--tree. In section \ref{sec_abelian},
starting with a free action of a group $\Gamma$ on an $\bbR^n$--tree $T$, we study the action on the $\bbR$--tree
$\ol T$ obtained by identifying points at infinitesimal distance, and we deduce a weaker version of the d\'evissage
Theorem where we obtain a graph of groups over (maybe non-finitely generated) abelian groups.
Section \ref{sec_flawless} contains the core of the argument: starting with a free action of $\Gamma$ on an $\bbR^n$--tree $T$,
we build a free action on an $\bbR^n$--tree $T'$ such that the $\bbR$--tree $\olT'$ has \emph{cyclic} arc stabilizers.
The d\'evissage theorem and its corollaries will then follow immediately, as shown in section \ref{sec_devissage}.

The author thanks the referee for useful presentation suggestions.

\section{Preliminaries}\label{sec_prelim}

\subsection{Marked groups and limit groups}\label{sec_prelim_limit}
Sela introduced limit groups in \cite{Sela_diophantine1}.
For background about Sela's limit groups, see also \cite{CG_compactifying} or \cite{Pau_theorie}.

A \emph{marked group} $(G,S)$ is a finitely generated group $G$ together with a finite ordered
generating family $S=(s_1,\dots,s_n)$. Note that repetitions may occur in $S$, and some generators $s_i$
may be the trivial element of $G$.
Consider two groups $G$ and $G'$ together with some markings of the same cardinality $S=(s_1,\dots,s_n)$
and $S'=(s'_1,\dots,s'_n)$.
A \emph{morphism of marked groups} $h\co (G,S)\ra(G',S')$ is a homomorphism 
$h\co G\ra G'$ sending $s_i$ on $s'_i$ for all $i\in\{1,\dots,n\}$. Note that
there is at most one morphism between two marked groups, and that all morphisms are epimorphisms.

A \emph{relation} in $(G,S)$ is an element of the kernel of the natural morphism $F_n\ra G$
sending $a_i$ to $s_i$ where $F_n$ is the free group with basis $(a_1,\dots,a_n)$.
Note that two marked group are isomorphic if and only if they have the same set of relations.

Given any fixed $n$, we define $\calg_n$ to be the set of isomorphism classes 
of marked groups. It is naturally endowed with the topology 
such that the sets $N_R(G,S)$ defined below form a neighbourhood basis of $(G,S)$.
For each $(G,S)\in\calg_n$ and each $R>0$, $N_R(G,S)$ is the set of marked groups $(G',S')\in\calg_n$
such that $(G,S)$ and $(G',S')$ have exactly the same relations of length at most $R$.
For this topology, $\calg_n$ is a Hausdorff, compact, totally disconnected space.

\begin{dfn}\rm
  A \emph{limit group} $(G,S)\in\calg_n$ is a marked group which is a limit of markings of free groups in $\calg_n$.
\end{dfn}

Actually, being a limit group does not depend on the choice of the generating set. 
Moreover, limit groups  have several equivalent characterizations: a finitely generated group is a limit group
if and only if it is fully residually free, if and only if it has the same universal theory as a free group,
if and only if it is a subgroup of a non-standard free group \cite{Remeslennikov_exist_siberian,CG_compactifying}.
We won't need those characterizations in this paper.

\subsection{$\Lambda$--trees}
For background on $\Lambda$--trees, see \cite{Bass_non-archimedean,Chi_book}.

\paragraph{Totally ordered abelian groups}

A totally ordered abelian group $\Lambda$ is an abelian group with a total ordering
such that for all $x,y,z\in \Lambda$, $x\leq y \Rightarrow x+z\leq y+z$.
Our favorite example will be $\bbR^n$, with the lexicographic ordering.
In all this paper, $\bbR^n$ will always be endowed with its lexicographic ordering.
To fix notations, we use the \emph{little endian} convention: 
the leftmost factor will have the greatest weight. More precisely, 
if $\Lambda_1$ and $\Lambda_2$ are totally ordered abelian groups,
the lexicographic ordering on $\Lambda_1\oplus\Lambda_2$ is
defined by $(x_1,x_2)\leq (y_1,y_2)$ if $x_1<y_1$ or ($x_1=y_1$ and $x_2\leq y_2$).

A morphism $\phi\co \Lambda\ra\Lambda'$ between two totally-ordered
abelian groups is a non-decreasing group morphism.
Given $a,b\in\Lambda$, the subset $[a,b]=\{x\in\Lambda \,|\, a\leq x\leq b\}$ is called the \emph{segment} between $a$ and $b$.
A subset $E\subset\Lambda$ is \emph{convex} if for all $a,b\in E$, $[a,b]\subset E$.
 The kernel of a morphism is a convex subgroup,
and if $\Lambda_0\subset \Lambda$ is a convex subgroup, then $\Lambda/\Lambda_0$ has a natural
structure of totally ordered abelian group. By \emph{proper} convex subgroup of $\Lambda$,
we mean a convex subgroup strictly contained in $\Lambda$. 

The set of convex subgroups of $\Lambda$ is totally ordered by inclusion.
The height of $\Lambda$ is the (maybe infinite) number of proper convex subgroups of $\Lambda$.
Thus, the height of $\bbR^n$ is $n$. $\Lambda$ is \emph{archimedean} if its height is at most 1.
It is well known that a totally ordered abelian group is archimedean if and only if it is isomorphic
to a subgroup of $\bbR$ (see for instance \cite[Theorem 1.1.2]{Chi_book})

If $\Lambda_0\subset \Lambda$ is a convex subgroup, then any element  $\lambda_0\in\Lambda_0$ may be thought as
infinitely small compared to an element $\lambda\in\Lambda\setminus\Lambda_0$ since for all $n\in\bbN$,
$n\lambda_0\leq \lambda$. Therefore, we will say that an element in $\bbR^n$ is \emph{infinitesimal}
if it lies in the maximal proper convex subgroup of $\bbR^n$, which we casually denote by $\bbR^{n-1}$.
Similarly, for $p\leq n$, we will identify $\bbR^p$ with the corresponding convex subgroup of $\bbR^n$.
The \emph{magnitude} of an element $\lambda\in\bbR^n$ is the smallest $p$ such that $\lambda\in \bbR^p$.
Thus $\lambda\in\bbR^n$ is infinitesimal if and only if its magnitude is at most $n-1$.

Given a totally ordered abelian group $\Lambda$, $\Lambda\tensor\bbQ$ has a natural structure of
a totally ordered abelian group by letting $\frac{\lambda}{n}\leq\frac{\lambda'}{n'}$ if and only if
$n'\lambda\leq n\lambda'$.

\paragraph{$\Lambda$--metric spaces and $\Lambda$--trees}

A \emph{$\Lambda$--metric space} $(E,d)$ is a set $E$ endowed with a map $d\co E\times E\ra \Lambda_{\geq 0}$
satisfying the three usual axioms of a metric: separation, symmetry and triangle inequality.
The set $\Lambda$ itself is a $\Lambda$--metric space for the metric $d(a,b)=|a-b|=\max(a-b,b-a)\in\Lambda$.
A \emph{geodesic segment} in $E$ is an isometric map from a segment $[a,b]\subset\Lambda$ to a subset of $E$.
A $\Lambda$--metric space is \emph{geodesic} if any two points are joined by a geodesic segment.
We will denote by $[x,y]$ a geodesic segment between two points in $E$ (which, in this generality, might be non-unique).

Note that even in a set $\Lambda$ like $\bbR^n$, the upper bound is not always defined so
one cannot easily define a $\Lambda$--valued diameter (see however \cite[p.99]{Chi_book} for
a notion of diameter as a interval in $\Lambda$).
Nevertheless, we will say that a subset $F$ of a $\bbR^n$--metric space $E$
is \emph{infinitesimal} if the distance between any two points of $F$ is infinitesimal.
Similarly, we define the \emph{magnitude} of $F$ as the smallest $p\leq n$ such that
the distance between any two points of $F$ has magnitude at most $p$.

We give two equivalent definitions of a $\Lambda$--tree. The equivalence is 
 proved for instance in \cite[Lemma 2.4.3,~p.71]{Chi_book}.
\begin{dfn}\rm
  A $\Lambda$--tree $T$ is a geodesic $\Lambda$--metric space such that
  \begin{itemize*}
    \item $T$ is $0$--hyperbolic in the following sense:
$$\forall x,y,u,v\in T,\  d(x,y)+d(u,v)\leq \max\{d(x,u)+d(y,v),d(x,v)+d(y,u)\}$$
\item $\forall x,y,z\in T,\ d(x,y)+d(y,z)-d(x,z) \in 2\Lambda$
  \end{itemize*}

Equivalently, a geodesic $\Lambda$--metric space is a $\Lambda$--tree if
\begin{itemize*}
\item the intersection of any two geodesic segments sharing a common endpoint is a geodesic segment
\item if two geodesic segments intersect in a single point, then their union is a geodesic segment.
\end{itemize*}
\end{dfn}

\begin{rem*}
In the first definition, the second condition is automatic if $2\Lambda=\Lambda$, which is the case for $\Lambda=\bbR^n$.

It follows from the definition that there is a unique geodesic joining a given pair of points in a $\Lambda$--tree.
\end{rem*}

Clearly, $\Lambda$ itself is $\Lambda$--tree.
Another simple example of a $\Lambda$--tree is the vertex set $V(S)$ of a simplicial tree $S$: 
$V(S)$ endowed with the combinatorial distance is a $\bbZ$--tree.

\subsection{Killing infinitesimals and extension of scalars}

The following two operations are usually known as the \emph{base change functor}.
\paragraph{Killing infinitesimals} Consider $\Lambda_0\subset \Lambda$
a convex subgroup (a set of infinitesimals), and let $\ol\Lambda=\Lambda/\Lambda_0$. If $\Lambda=\bbR^n$, we will usually take $\Lambda_0=\bbR^{n-1}$,
so that $\ol\Lambda\simeq\bbR$. Consider a $\Lambda$--metric space $E$. Then the relation $\sim$ defined by
$x\sim y \Leftrightarrow d(x,y)\in \Lambda_0$ is an equivalence relation on $E$, and the $\Lambda$--metric on $E$
provides a natural $\ol\Lambda$--metric on $E/\sim$. 
We say that $\ol E=E/\sim$ is obtained from $E$ by \emph{killing infinitesimals}.
Clearly, if $T$ is a $\Lambda$--tree, then $\ol T$ is a $\ol\Lambda$--tree.
Thus, killing infinitesimals in an $\bbR^n$--tree $T$ provides an $\bbR$--tree $\ol T$.
By extension, we will often denote $\bbR$--trees with a bar.

\paragraph{Extension of scalars}
Consider a $\Lambda$--tree $T$, and an embedding $\Lambda\hookrightarrow \Tilde\Lambda$
(for example, one may think of $\bbZ\subset\bbR$).
Then $T$ may be viewed as a $\Tilde\Lambda$--metric space, but it is not 
a $\Tilde\Lambda$--tree if $\Lambda$ is not convex in $\Tilde\Lambda$: as a matter of fact, $T$ is not
geodesic as a $\Tilde\Lambda$--metric space (there are holes in the geodesics).
However, there is a natural way to \emph{fill the holes}:
\begin{prop}[Extension of scalars, see {\cite[Theorem 4.7, p.75]{Chi_book}}]
  There exists a $\Tilde\Lambda$--tree $\Tilde T$ and an isometric embedding $T\hookrightarrow \Tilde T$
which is canonical in the following sense:
 if $T'$ is another $\Tilde\Lambda$--tree with an isometric embedding $T\hookrightarrow T'$,
then there is a unique $\Tilde\Lambda$--isometric embedding $\Tilde T\ra T'$ commuting with the embeddings of $T$
in $\Tilde T$ and $T'$.
\end{prop}

For example, take $T$ to be the $\bbZ$--tree corresponding corresponding to the set of vertices of a simplicial tree $S$.
Then the embedding $\bbZ\subset \bbR$ gives an $\bbR$--tree $\Tilde T$ which is isometric to the geometric
realization of $S$.

\begin{rem*}
  The proposition also holds if one only assumes that $T$ is $0$--hyperbolic. In this case, taking $\Tilde\Lambda=\Lambda$,
one gets a natural $\Lambda$--tree containing $T$.
\end{rem*}

\subsection{Subtrees}

A \emph{subtree} $Y$ of a $\Lambda$--tree $T$ is a convex subset of $T$, \ie such that for all $x,y\in Y$,
$[x,y]\subset Y$. A subtree is \emph{non-degenerate} if it contains at least two points.
One could think of endowing $\Lambda$, and $T$, with the order topology.
However, this is usually not adapted. For instance: $\bbR^n$ is not connected
with respect to this topology for $n>1$. 
This is why we need a special definition of a \emph{closed} subtree.
The definition coincides with the topological definition for $\bbR$--trees.

\begin{dfn}[Closed subtree]\rm
  A subtree $Y\subset T$ is a \emph{closed subtree} if the intersection of $Y$ with a segment of $T$
is either empty or a segment of $T$.
\end{dfn}

There is a natural projection on a closed subtree. Consider a base point $y_0\in Y$.
Then for any point $x\in T$, there is a unique point $p\in Y$ such that $[y_0,x]\cap Y=[y_0,p]$.
One easily checks that $p$ does not depend on the choice of the base point $y_0$
($[p,x]$ is the \emph{bridge} between $x$ and $Y$, see \cite{Chi_book}).
The point $p$ is called the \emph{projection} of $x$ on $Y$. 

\begin{rem*}
The existence of a projection is actually equivalent to the fact that the subtree $Y$ is closed.
Be aware that a non-trivial proper convex subgroup of $\Lambda$ is never closed in $\Lambda$.
In particular, the intersection of infinitely many closed subtrees may fail to be closed.
\end{rem*}

A \emph{linear} subtree of $T$ is a subtree in which any three points are contained in a segment.
It is an easy exercise to prove that a maximal linear subtree of $T$ is closed in $T$.
Finally, any linear subtree $L\subset T$ is isometric to a convex subset of $\Lambda$
and any two isometries $L\ra \Lambda$ differ by an isometry of $\Lambda$.

\subsection{Isometries}
An isometry $g$ of a $\Lambda$--tree $T$ can be of one of the following exclusive types:
\begin{itemize*}
  \item elliptic: $g$ has a fix point in $T$
  \item inversion: $g$ has no fix point, but $g^2$ does
  \item hyperbolic: otherwise. 
\end{itemize*}
In all cases, the set $A_g=\{x\in T\,|\, [g\m x,x]\cap[x,g.x]=\{x\}$
is called the \emph{characteristic set} of $g$.

If $g$ is elliptic, $A_g$ is the set of fix points of $g$ which is a closed subtree of $T$.
Moreover, for all $x\in T$, the midpoint of $[x,g.x]$ exists and lies in $A_g$.

If $g$ is an inversion, then $A_g=\es$. Actually, for any $x\in T$, $d(x,g.x)\notin 2\Lambda$
so $[x,g.x]$ has no midpoint in $T$. In particular, if $2\Lambda=\Lambda$
(which occurs for instance if $\Lambda=\bbR^n$), inversions don't exist.
Moreover, one can perform the analog of barycentric subdivision for simplicial trees
to get rid of inversions: consider $\Tilde\Lambda=\frac12 \Lambda$,
and let $\Tilde T$ be the $\Tilde\Lambda$--tree obtained by the extension of scalars $\Lambda\subset\Tilde\Lambda$.
Then the natural extension of $g$ to $\Tilde T$ fixes a unique point in $\Tilde T$ (in particular, $g$
is elliptic in $\Tilde T$).
If $g$ is elliptic or is an inversion, its translation length $l_T(g)$ is defined to be $0$.

If $g$ is hyperbolic, then the set $A_g$
is non-empty, and is a maximal linear subtree of $T$, and is thus closed in $T$. 
It is called \emph{the axis} of $g$. 
Moreover, the restriction of $g$ to $A_g$
is conjugate to the action of a translation $\tau\co a\mapsto a+l_T(g)$ on a $\tau$--invariant
convex subset of $\Lambda$ for some positive $l_T(g)\in\Lambda$. The \emph{translation length}
 of $g$ is the element $l_T(g)\in\Lambda_{>0}$.
If $p$ is the projection of $x$ on $A_g$, then for $k\neq 0$, $d(x,g^k.x)=2d(x,p)+|k|l_T(g)$.

Note that it may happen that $A_g$ is not isometric to $\Lambda$.
For instance, if $\Lambda=\bbR^2$, the axis of an element $g$ with infinitesimal translation length
can be of the form $I\times \bbR$ where $I$ is any non-empty interval in $\bbR$ which can be
open, semi-open or closed.

If $g$ is hyperbolic, then for all $x\in T$, the projection of $x$ on $A_g$ is the projection of $x$ on $[g\m.x,g.x]$.
In particular, if the midpoint of $[x,g.x]$ exists, then it lies in $A_g$. It also follows that if
$g$ is hyperbolic and if $g\m.x,x,g.x$ are aligned (in any order) then they lie on the axis of $g$.

If an abelian group $A$ acts by isometries on $\Lambda$--tree $T$ and contains
a hyperbolic element $g$, then all the hyperbolic elements of $A$
have the same axis $l$, $A$ contains no inversion, and all elliptic elements
fix $l$. We say that $l$ is the \emph{axis} of the abelian group $A$.
The axis of $A$ can be characterized as the only \emph{closed}
$A$--invariant linear subtree of $T$, or as the only \emph{maximal}
$A$--invariant linear subtree of $T$.

Finally if a group $\Gamma$ has a non-abelian action without inversion on a $\Lambda$--tree $T$, 
then there is a unique minimal non-empty $\Gamma$--invariant subtree $\Tmin(\Gamma)$ \cite[Theorem  3.4.1]{Chi_introduction}.
This applies if $\Gamma$ is non-abelian and acts freely on $T$.
Note however that $\Tmin(\Gamma)$ is not a closed subtree of $T$ in general (even if $\Lambda=\bbR$);
if $\Lambda$ is not archimedean, $\Tmin(\Gamma)$ may also fail to coincide with the union of
the translation axes of $\Gamma$;
finally, there is no unique such minimal subtree for the action of $\Gamma=\{0\}\times \bbZ$ on $\bbR^2$.

\subsection{Elementary properties of groups acting freely on $\Lambda$--trees}

We now recall some elementary properties of groups acting freely (without inversion)
on $\Lambda$--trees. They are proved for instance in \cite{Chi_book}.
\begin{lem}\label{lem_elementary}
  Let $\Gamma$ be a group acting freely without inversion on a $\Lambda$--tree.
Then 
\begin{enumerate*}
\item $\Gamma$ is torsion free;
\item two elements $g,h\in\Gamma$ commute if and only if they have the same axis. If they don't commute,
the intersection of their axes is either empty or a segment (\cite{Chi_book}, proof of Lemma  5.1.2 p.218 and Rk p.111);
\item maximal abelian subgroups of $\Gamma$ are malnormal
(property CSA) and $\Gamma$ is commutative transitive:
the relation of commutation on $\Gamma\setminus \{1\}$ is transitive \cite[Lemma  5.1.2 p.218]{Chi_book}.
\end{enumerate*}
\end{lem}

\begin{rem*}
Property CSA implies that $\Gamma$ is commutative transitive.

A result known as \emph{Harrison Theorem}, proved by Harrison for $\bbR$--trees 
and by Chiswell and Urbanski--Zamboni for general $\Lambda$--trees, 
says that a $2$--generated group acting freely without inversion on a $\Lambda$--tree is either a free group or a free abelian group.
(see \cite{Chi_harrison,Urbanski-Zamboni_Free,Harrison_real}). We won't use this result in this paper.
\end{rem*}

\section{A limit group acts freely on an $\bbR^n$--tree}\label{sec_remeslennikov}

The goal of this section is to prove that limit groups act freely on $\bbR^n$--trees.
This is an adaptation of an argument by Remeslennikov concerning fully
residually free groups (\cite{Remeslennikov_exist_ukrainian}, see also \cite[Theorem 5.5.10~p.246]{Chi_book}).
Note that it is claimed in \cite{Remeslennikov_exist_ukrainian}
that finitely generated fully residually free groups act freely on a $\Lambda$--tree where
$\Lambda$ is a \emph{finitely generated} ordered abelian group.
However, the proof is not completely correct since it relies on a misquoted result about
valuations (Theorem 3 in \cite{Remeslennikov_exist_ukrainian}) 
to which there are known counterexamples 
(for any subgroup $\Lambda\subset\bbQ$, there is valuation on $\bbQ(X,Y)$, extending the trivial valuation on $\bbQ$,
 whose value group is $\Lambda$ \cite[ch.VI, \para 15,~ex.3,4]{Zariski-Samuel_commutative2} or \cite[Theorem 1.1]{Kuhlmann_value}). 
Nevertheless, Remeslennikov's argument proves
 the following weaker statement: a finitely generated 
fully residually free group acts freely on a $\Lambda$--tree where $\Lambda$ has finite $\bbQ$--rank, 
\ie $\Lambda\tensor\bbQ$ is finite dimensional over $\bbQ$.

The fact that a limit group acts freely  on an $\bbR^n$ tree
 will be deduced from a more general result about group acting freely on Bruhat--Tits trees.
But we first state a simpler result in this spirit (see also \cite{GaSp_does,GaSp_every}).
Remember that $\calg_n$ denotes the space of groups marked by a generating family of cardinality $n$.

\begin{prop}[Acting freely on $\Lambda$--trees is closed]
Let $\calt_n\subset\calg_n$ be the set of marked groups having a free action without inversion on some $\Lambda$--tree
($\Lambda$ may vary with the group).

Then $\calt_n$ is closed in $\calg_n$.
\end{prop}

We won't give the proof of this result since this proposition is not sufficient for us 
as it does not give any control over $\Lambda$.
This is why we rather prove the following more technical result.%
\footnote{The proof is actually  very similar
to the proof of the more technical result: instead of taking ultraproducts of valuated fields, take an ultraproduct
of trees to get a free action without inversion on a $\Lambda^*$--tree
(see also \cite[p.239]{Chi_book} where the behavior $\Lambda$--trees under ultrapowers is studied
in terms of Lyndon length functions).}

For general information of the action of $SL_2(K)$ on its Bruhat--Tits $\Lambda$--tree $BT_K$ where $K$ a field,
and $v\co K\ra\Lambda\cup\{\infty\}$ is a valuation,
see for instance \cite[\para 4.3, p.144]{Chi_book}.
 Essentially, we will only use the existence of the Bruhat--Tits $\Lambda$--tree and
the formula for the translation length of a matrix $m\in SL_2(K)$:
$l_{BT_K}(m)=\max\{-2v(\Tr(m)),0\}$. Also note that the action of $SL_2(K)$ on its
Bruhat--Tits tree has no inversion (however, there may be inversions in $GL_2(K)$).

\begin{dfn}[Action on a Bruhat--Tits tree]\rm
  By an action of $\Gamma$ on a Bruhat--Tits tree, we mean an action of $\Gamma$ 
on the Bruhat--Tits $\Lambda$--tree for $SL_2(K)$ induced by a morphism $j\co \Gamma\ra SL_2(K)$ 
where $K$ is a valuated field with values in $\Lambda$.
\end{dfn}

\begin{SauveCompteurs}{thm_closed}%
\thmclosed
\end{SauveCompteurs}

\begin{SauveCompteurs}{Remeslennikov}%
\corRemeslennikov
\end{SauveCompteurs}

\begin{proof}[Proof of the corollary]
This follows from the theorem above since a free group acts freely on a Bruhat--Tits tree.
\end{proof}

\begin{proof}[Proof of the Theorem]
We first prove that $\calbt$ is closed.
Let $(\Gamma_i,S_i)\in\calbt$ be a sequence of marked groups converging to $(\Gamma,S)$.
For each index $i$, consider a field $K_i$ and a valuation $v_i\co K_i\ra \Lambda_i\cup\{\infty\}$
and an embedding $j_i\co \Gamma_i\ra SL_2(K_i)$ such that $j_i(\Gamma_i)$ acts freely without inversions
on the corresponding Bruhat--Tits tree $BT_i$.

  Consider $\omega$ an ultrafilter on $\bbN$, \ie 
 a finitely additive measure of total mass $1$ (a \emph{mean}), defined on all subsets of $\bbN$,
and with values in $\{0,1\}$, and assume that this ultrafilter is non-principal, \ie 
that the mass of finite subsets is zero. 
Say that a property $P(k)$ depending on $k\in\bbN$ is true \emph{$\omega$--almost everywhere}
if $\omega(\{k\in\bbN|P(k)\})=1$. Note that a property which is not true almost everywhere 
is false almost everywhere.
Given a sequence of sets $(E_i)_{i\in\bbN}$, the ultraproduct $E^*$ of $(E_i)$
is the quotient $(\prod_{i\in \bbN} E_i)/\eqv_\omega$ 
where  $\eqv_\omega$ is the natural equivalence relation on $\prod_{i\in \bbN} E_i$
defined by equality $\omega$--almost everywhere.

Consider $K^*$ the ultraproduct of the fields $K_i$,
$\Gamma^*$ the ultraproduct of the groups $\Gamma_i$, and
$\Lambda^*$  the ultraproduct of the totally ordered abelian groups $\Lambda_i$.
As a warmup, we prove the easy fact that the natural ring structure on $K^*$ makes it a field:
if $k^*=(k_i)_{i\in\bbN}\neq 0$ in $K^*$, then for almost all $i\in\bbN$,
$k_i\neq 0$, and $1/k_i$ is defined for almost every index $i$, and defines an inverse $(1/k_i)_{i\in\bbN}$ for $k^*$ in $K^*$.

Similarly, $\Gamma^*$ is a group, and $\Lambda^*$ a totally ordered abelian group (for the total 
order $(x_i)_{i\in\bbN}\leq (y_i)_{i\in\bbN}$ if and only if $x_i\leq y_i$ almost everywhere).
Now consider the map $v^*\co K^*\ra\Lambda^*\cup\infty$ defined by $v^*((k_i)_{i\in\bbN})=(v_i(k_i))_{i\in\bbN}$,
and the map $j^*\co \Gamma^*\ra SL_2(K^*)$ defined by $j^*((g_i)_{i\in\bbN})=(j_i(g_i))_{i\in\bbN}$.
Then $v^*$ is a valuation on $K^*$, and $j^*$ a monomorphism of groups.
We denote by $BT^*$ the Bruhat--Tits tree of $SL_2(K^*)$.(\footnote{It may also be checked that 
$BT^*$ is actually the ultraproduct of the $\Lambda_i$--trees $BT_i$.})

Now, given a field $K$ with a valuation $v\co K\ra \Lambda\cup\{\infty\}$, a subgroup $H\subset SL_2(K)$
acts freely without inversions on the corresponding Bruhat--Tits tree $BT$ if and only if
the translation length of any element $h\in H\setminus\{1\}$ is non-zero.
But the translation length of a matrix $m\in SL_2(K)$ can be computed in terms of the valuation of its trace by
the formula $l_{BT}(m)=\max\{0,-2v(\Tr(m))\}$, so the freeness (without inversion) of the action translates into
$v(\Tr(h))<0$ for all $h\in H\setminus\{1\}$ \cite[Lemma 4.3.5~p.148]{Chi_book}.
Therefore, since for all $i$ and all $g_i\in \Gamma_i\setminus\{1\}$,
$\Tr(j_i(g_i))$ has negative valuation,
all the elements $g^*\in \Gamma^*\setminus\{1\}$ satisfy  $v^*(\Tr(j^*(g^*)))<0$,
which means that $\Gamma^*$ acts freely without inversion on $BT^*$.

Finally, there remains to check that the marked group $(\Gamma,S)$ embeds into $\Gamma^*$ (see for instance \cite{CG_compactifying}).
We use the notation $S=(s_1,\dots,s_n)$ and $S_i=(s_1^{(i)},\dots,s_n^{(i)})$.
Consider the family $S^*=(s_1^*,\dots,s_n^*)$ of elements of $\Gamma^*$ defined by
$s_1^*=(s_1^{(i)})_{i\in\bbN}, \dots,s_n^*=(s_n^{(i)})_{i\in\bbN}$.
The definition of the convergence of marked groups says that
if an $S$--word represents the trivial element (resp.\ a non-trivial element) in $\Gamma$, then
for $i$ sufficiently large, the corresponding $S_i$--word is trivial (resp.\ non-trivial) in $\Gamma_i$.
Since $\omega$ is non-principal, this implies that the corresponding $S^*$--word is trivial (resp.\ non-trivial).
This means that the map sending $(s_1,\dots,s_n)$ to $(s_1^*,\dots,s_n^*)$ extends to an isomorphism
between $\Gamma$ and $\langle S^* \rangle\subset \Gamma^*$.
Therefore, $(\Gamma,S)\in\calbt$, so $\calbt$ is closed.

We now prove the fact that any group $(\Gamma,S)$ in $\calbt$ acts freely on a $\Lambda$--tree where
$\Lambda\tensor\bbQ$ has dimension at most $3n+1$. So consider an embedding $j\co \Gamma\ra SL_2(K)$ where 
$K$ has a valuation $v\co K\ra\Lambda\cup\{\infty\}$ such that the induced action of $\Gamma$ on the 
Bruhat--Tits tree for $SL_2(K)$ is free without inversion.
Consider the subfield $L\subset K$ generated by the $4n$ coefficient of the matrices
$j(s_1),\dots,j(s_n)$. Since the matrices have determinant $1$, $L$ can be written
as $L=k_0(x_1,\dots,x_{3n})$ where $k_0$ is the prime subfield of $K$. Let $\Lambda_L=v(L\setminus\{0\})$ be the value group of $L$.
Since $\Gamma$ embeds in $SL_2(L)$, $\Gamma$ acts freely on the corresponding Bruhat--Tits $\Lambda_L$--tree.
We now quote a result about valuations which implies that $\Lambda_L$ has finite $\bbQ$--rank.

\begin{thm}{\rm\cite[Corollary 1 in VI.10.3]{Bourbaki_valuations}}\label{thm_Bourbaki}\qua
  Let $L=L_0(x_1,\dots,x_p)$ be a finitely generated extension of $L_0$, and $v\co L\ra \Lambda\cup\{\infty\}$ a valuation.
Denote by $\Lambda_L=v(L\setminus\{0\})$ (resp.\ $\Lambda_0=v(L_0\setminus\{0\})$) the corresponding value group.
Then the $\bbQ$--vector space $(\Lambda_L\tensor \bbQ)/(\Lambda_0\tensor \bbQ)$ has dimension at most $p$.
\end{thm}

Taking $L_0=k_0$, one gets that $\Lambda_L$ has $\bbQ$--rank at most $3n+1$ since 
$\Lambda_0$ is either trivial or isomorphic to $\bbZ$.

Using the extension of scalars (base change functor), there remains to prove that 
if a totally ordered group $\Lambda$ has finite $\bbQ$--rank, then it is isomorphic to a subgroup of
$\bbR^n$.
\end{proof}

\begin{lem}
Consider $\Lambda$ a totally ordered group of $\bbQ$--rank $p$.
Then $\Lambda$ is isomorphic (as an ordered group) to a subgroup of $\bbR^p$ with its lexicographic ordering.
\end{lem}

\begin{rem*}
  However, $\Lambda\tensor\bbQ$  is usually not isomorphic to $\bbQ^p$ with its lexicographic ordering
as shows an embedding of $\bbQ^2$ into $\bbR$.
\end{rem*}

\begin{proof}
We first check that the height of $\Lambda$ is at most $p$ (see \cite[prop 3 in VI.10.2]{Bourbaki_valuations}).
  First, $\Lambda$ embeds into $\Lambda\tensor\bbQ$, so we may replace $\Lambda$ by $\Lambda\tensor\bbQ$
and assume that $\Lambda$ is a totally ordered $\bbQ$--vector space of dimension $p$.
Any convex subgroup $\Lambda_0\subset\Lambda$ is a $\bbQ$ vector subspace in $\Lambda$
since if $0\leq x\in\Lambda_0$, for all $k\in\bbN\setminus\{0\}$,
$\frac{1}{k}x \in \Lambda_0$ since $0\leq \frac{1}{k}x \leq x$.
Now the height of $\Lambda$ is at most $p$ since a chain of convex subgroups 
$\Lambda_0\subsetneqq \Lambda_1 \subsetneqq \dots \subsetneqq \Lambda_i$ 
is a chain of vector subspaces.

We now prove by induction that a totally ordered group $\Lambda$ of finite height $q$ embeds as an ordered subgroup 
of $\bbR^q$ with its lexicographic ordering. Once again, one can replace $\Lambda$ by $\Lambda\tensor\bbQ$ without 
loss of generality. We argue by induction on the height.
If $\Lambda$ has height $1$, \ie if $\Lambda$ is archimedean, then $\Lambda$ embeds in $\bbR$ 
(see for instance \cite[Theorem 1.1.2]{Chi_book}).
Now consider $\Lambda_0\subset\Lambda$ the maximal proper convex subgroup of $\Lambda$, and let 
$\ol\Lambda=\Lambda/\Lambda_0$.
Since $\Lambda,\Lambda_0$ and $\ol\Lambda$ are $\bbQ$ vector spaces, one has algebraically
$\Lambda= \ol \Lambda \oplus \Lambda_0$.

The fact that $\Lambda_0$ is convex in $\Lambda$ implies
that the ordering on $\Lambda$ corresponds to the lexicographic ordering on $\Lambda\oplus \Lambda_0$
\cite[lemma 2 in VI.10.2]{Bourbaki_valuations}.
Indeed, one first easily checks that any section $j\co \ol\Lambda\ra\Lambda$ 
is increasing.
Now let's prove that the
isomorphism $f\co \ol\Lambda\times \Lambda_0\ra\Lambda$ defined by  $f(\ol x,x_0)=j(\ol x)+x_0$ is increasing
for the lexicographic ordering on $\ol\Lambda\times \Lambda_0$.
So assume that $(\ol x,x_0)\geq 0$. If $\ol x=0$, then $f(\ol x,x_0)=x_0\geq 0$.
If $\ol x>0$, then $f(\ol x,x_0)=j(\ol x)+x_0>0$ since
otherwise, one would have $0\leq j(\ol x)\leq -x_0$, hence $j(\ol x)\in \Lambda_0$ by convexity, a contradiction.

Finally, by induction hypothesis, $\Lambda_0$ embeds as an ordered subgroup of $\bbR^{q-1}$ and $\ol\Lambda$
embeds as an ordered subgroup of $\bbR$, so $\Lambda$ embeds as an ordered subgroup of $\bbR^q$.
\end{proof}

\section{Gluing $\Lambda$--trees}\label{sec_gluing}
The goal of this section is to define graph of actions on $\Lambda$--trees which show how to glue actions on $\Lambda$--trees
along closed subtrees to get another action on a $\Lambda$--tree, and to give a criterion for the resulting action to be free.
We will finally study more specifically gluings of $\bbR$--trees along points, and show that a decomposition
of an $\bbR$--tree $T$ into a graph of actions on $\bbR$--trees above points correspond to a \emph{transverse} covering
of $T$ by closed subtrees.

\subsection{Gluing $\Lambda$--trees along points}
Here, we recall that one can glue $\Lambda$ trees together along a point
to get a new $\Lambda$--tree (see \cite[Lemma 2.1.13]{Chi_book}).

\begin{lem}{\rm\cite[Lemma 2.1.13]{Chi_book}}\label{lem_points}\qua
  Let $(Y,d)$ be a $\Lambda$--tree, and $(Y_i,d_i)_{i\in I}$ be a family of $\Lambda$--trees.
Assume that $Y_i\cap Y=\{x_i\}$ and that for all $i,j\in I$, $Y_i\cap Y_j=\{x_i\}\cap\{x_j\}$. 
Let $X=(\bigcup_{i\in I}Y_i)\cup Y$
and let $\ol{d}\co X\times X\ra\Lambda$ defined by: $\ol{d}_{|Y\times Y}=d$;
$\ol{d}_{|Y_i\times Y_i}=d_i$; for $x\in Y_i,y\in Y$ $\ol{d}(x,y)=d(x,x_i)+d(x_i,y)$;
for $x\in Y_i,y\in Y_j$ $\ol{d}(x,y)=d(x,x_i)+d(x_i,x_j)+d(x_j,y)$.

Then $(X,\ol{d})$ is a $\Lambda$--tree.
\end{lem}

\subsection{Gluing two trees along a closed subtree}
The following gluing construction will be used for gluing trees along 
maximal linear subtrees.

Assume that we are given two $\Lambda$--trees $(Y_1,d_1),(Y_2,d_2)$, 
two \emph{closed} subtrees $\delta_1\subset Y_1$ and $\delta_2\subset Y_2$,
and an isometric map $\phi\co \lambda_1\onto\lambda_2$.
By definition of a closed subtree, we have two orthogonal projections
$p_{\lambda_i}\co Y_i\ra \lambda_i$ for $i\in\{1,2\}$.

Let $X=Y_1 \dunion Y_2$, and let $\sim$ be the equivalence relation on $X$
generated by $x\sim \phi(x)$ for all $x\in \lambda_1$.
The set $(Y_1\cup_\phi Y_2):=X/\sim$ is now endowed with the following metric
which extends $d_i$ on $Y_i$:
if $x\in Y_1$ and $y\in Y_2$, we set 
\begin{eqnarray}
  d(x,y)&:=&d_1(x,p_{\lambda_1}(x))+d_2(\phi(p_{\lambda_1}(x)),p_{\lambda_2}(y))+d_2(y,p_{\lambda_2}(y))\label{eq_proj}\\
&=&d_1(x,p_{\lambda_1}(x))+d_2(\phi(p_{\lambda_1}(x)),y) \notag \\
&=&\min \{d_1(x,x_1)+d_2(\phi(x_1),y)\ | \ x_1\in\lambda_1\}\label{eq_min}
\end{eqnarray}

To prove the last equality, introduce $z_1=p_{\lambda_1}(x)$; then for any $x_1\in\lambda_1$,
\begin{eqnarray*}
d_1(x,x_1)+d_2(\phi(x_1),y)&=&d_1(x,z_1)+d_1(z_1,x_1)+d_2(\phi(x_1),y)\\
&=&d_1(x,z_1)+d_1(\phi(z_1),\phi(x_1))+d_2(\phi(x_1),y)\\
&\geq& d_1(x,z_1)+d_1(\phi(z_1),y)\\
&=&d(x,y)
\end{eqnarray*}

\begin{lem}\label{lem_glue2}
  With the definitions above, $(Y_1\cup_\phi Y_2,d)$ is a $\Lambda$--tree.
Moreover, any closed subtree of $Y_i$ is closed in $(Y_1\cup_\phi Y_2,d)$.
\end{lem}

\begin{proof}
Let $T=Y_1\cup_\phi Y_2$.
Then $T$ can be viewed as 
the tree $L=\lambda_1=\lambda_2$ on which are glued some subtrees of $Y_1,Y_2$ at some points.
More precisely, 
for $x\in \lambda_1$, let $A_x=(p_{\lambda_1})\m(x)$, 
and similarly, for $x\in \lambda_2$, let $B_x=(p_{\lambda_2})\m(x)$.
Then, because of the formula \eqref{eq_proj} for the metric, 
$T$ is isometric to the $\Lambda$--tree obtained by gluing
the trees $A_x$ and $B_x$ on $L$ along the point $x$ as in lemma \ref{lem_points}.
Therefore, by lemma \ref{lem_points},  $(T,d)$ is an $\Lambda$--tree.

Now let $Z\subset Y_1$ be a closed subtree. Let's prove that $Z$ is closed in $T$.
Consider $z\in Z$, $y\in T$, and let's prove that there exists $z_0\in Z$ such that
$[z,y]\cap Z=[z,z_0]$. If $y\in Y_1$, then one can take $z_0$ to be the projection of $y$ on $Z$
by hypothesis. If $y\in Y_2$, let $y_0$ be the projection of $y$ on $\lambda_2$,
and $z_0$ be the projection of $\phi\m(y_0)$ on $Z$. 
Of course, $[y,z_0]\subset Z$.
Now $[y_0,y]\setminus\{y_0\}$ does not meet $Z$ since it is contained in $T\setminus Y_1$,
and neither does $[z_0,y_0]\setminus\{z_0\}$.
Thus $[z,y]\cap Z=[z,z_0]$ and $Z$ is closed in $T$.
\end{proof}

\subsection{Equivariant gluing: graphs of actions on $\Lambda$--trees}

The combinatorics of the gluing will be given by a simplicial tree $S$,
 endowed with an action without inversion of a group $\Gamma$.
We denote by $V(S)$ and $E(S)$ the set of vertices and (oriented) edges of $S$, by $t(e)$ and $o(e)$
the origin and terminus of an (oriented) edge $e$, and by $\ol e$ the edge with opposite orientation as $e$.

A graph of actions on trees is usually defined as a graph of groups with some additional data like vertex trees.
Here, we rather use an equivariant definition at the level of the Bass--Serre tree.

\begin{dfn}[Graph of actions on $\Lambda$--trees]\rm
  Given a group $\Gamma$, a \emph{$\Gamma$--equiv\-a\-iant graph of actions on $\Lambda$--trees} 
is a triple $(S,(Y_v)_{v\in V(S)},(\phi_e)_{e\in E(S)})$
where
  \begin{itemize*}
  \item $S$ is a simplicial tree,
  \item for each vertex $v\in V(S)$, $Y_v$ is a $\Lambda$--tree (called \emph{vertex tree}), 
  \item for each edge $e\in E(S)$, $\phi_e\co \lambda_{\ol e}\onto \lambda_e$ is an isometry
between closed subtrees $\lambda_{\ol e}\subset Y_{o(e)}$ and $\lambda_e\subset Y_{t(e)}$
such that $\phi_{\ol e}=\phi_e\m$. We call the subtrees $\lambda_e$ the \emph{edge subtrees}.
  \end{itemize*}

This data is assumed to be $\Gamma$--\emph{equivariant} in the following sense:
\begin{itemize*}
\item $\Gamma$ acts on $S$ without inversion,
\item $\Gamma$ acts on $X=\dunion_{v\in V(S)} Y_v$ so 
  that the restriction of each element of $\Gamma$ to a vertex tree is an isometry,
\item the natural projection $\pi\co X\ra V(S)$ (sending a point in $Y_v$ to $v$)
is equivariant
\item the family of gluing maps is equivariant: 
for all $g\in\Gamma$, $\lambda_{g.e}=g.\lambda_e$, 
and $\phi_{g.e}=g\circ\phi_e\circ g\m$.
\end{itemize*}
\end{dfn}

We say that such a graph of actions is \emph{over} its edge subtrees.

\paragraph{The $\Lambda$--tree dual to a graph of actions}
Given $\calg$ a $\Gamma$--equivariant graph of actions on $\Lambda$--trees, 
we consider the smallest equivalence relation $\sim$ on $X=\dunion_{v\in V(S)} Y_v$
such that for all edge $e\in E(S)$ and $x\in \lambda_{\ol e}$, $x\sim \phi_e(x)$.
The $\Lambda$--tree dual to $\calg$ is the quotient space $T_\calg=X/\!\!\sim$.
To define the metric on $T_\calg$, one can alternatively say that $T_\calg$ is obtained
by gluing successively the vertex trees along the edge trees
according to lemma \ref{lem_glue2} in previous section.
Formula \eqref{eq_min} in previous section shows that the metric does not depend on the order in which the gluing are performed.
Indeed, an induction shows that the distance between $x\in Y_u$ and $y\in Y_v$
can be computed as follows:
let $e_1,\dots,e_n$ the edges of the path from $u$ to $v$ in $S$,
and $v_0=u,v_1,\dots,v_n=v$ the corresponding vertices
then
$$d(x,y)=\min\{ d_{Y_u}(x,x_1)+d_{Y_{v_1}}(\phi_{{e_1}}(x_1),x_2)+\dots+
d_{Y_{v_n}}(\phi_{{e_n}}(x_n),y)$$
where the minimum is taken over all $x_i\in \lambda_{\ol e_i}$.
By finitely many applications
of lemma \ref{lem_glue2}, one gets that the gluings corresponding to finite subtrees
of $S$ are $\Lambda$--trees. Now apply the fact that an increasing union of $\Lambda$--trees
is a $\Lambda$--tree to get that $T$ is a $\Lambda$--tree (see \cite[Lemma 2.1.14]{Chi_book}).
We thus get the following lemma:

\begin{dfn}[Tree dual to a graph of actions on $\Lambda$--trees]\label{dfn_goa}
\rm Consider $\calg=(S,(Y_v),(\phi_e))$ a $\Gamma$--equivariant graph of actions on $\Lambda$--trees.
The dual tree $T_\calg$ is the set $X/\!\!\sim$
endowed with the metric $d$ defined above.
It is a $\Lambda$--tree on which $\Gamma$ acts by isometries.

We say that a $\Lambda$--tree $T$ \emph{splits} as a graph of actions $\calg$ if
there is an equivariant isometry between $T$ and $T_\calg$.
\end{dfn}

\begin{rem*}
Consider an increasing union of trees $T_i$ such that $Y\subset T_0$ 
is a closed subtree of each $T_i$. Then $Y$ is closed in $\cup_i T_i$. Therefore,
using lemma \ref{lem_glue2}, one gets that a closed subtree of a vertex tree is closed in $T_\calg$.
In particular, vertex trees themselves are closed in $T_\calg$.
\end{rem*}

\subsection{Gluing free actions into free actions}
We next give a general criterion saying that an action obtained by gluing is free.
It is stated in terms of the equivalence relation $\sim$ on 
$X=\dunion_{v\in V(S)} Y_v$ defined above.
Each equivalence class has a natural structure of a connected graph:
elements of the equivalence class are vertices, put an oriented edge between two vertices
$x$ and $y$ if $y=\phi_e(x)$ for some edge $e\in E(S)$. 
Since this graph embeds into $S$ via the map
$\pi\co X\ra S$, so this graph is a simplicial tree. 
This graph structure defines a natural combinatorial metric on each equivalence class.

\begin{lem}[Criterion for a graph of free actions to be free] \label{lem_criterion_free}
Consider $\calg=(S,Y_v,\phi_e)$  a $\Gamma$--equivariant graph of actions on $\Lambda$--trees.
For each vertex $v\in V(S)$, denote by $\Gamma_v$ its stabilizer, and assume that
the action of $\Gamma_v$ on $Y_v$ is free. Assume furthermore that
 each equivalence class of $\sim$ has finite diameter (for the combinatorial metric).

Then the action of $\Gamma$ on $T_\calg$ is free.
\end{lem}

\begin{rem*}
  If the equivalence classes are allowed to have infinite diameter, then the resulting
action may fail to be free. The simplest example is as follows: take $S$ a simplicial 
tree with a free action of a non-trivial free group $\Gamma$, for instance $S$ is a line and $\Gamma=\bbZ$.
Define each vertex tree $Y_v$ as a point (having a free action of the trivial group), and $\phi_e$
the only map. Then $T_\calg$ is a point.

One can cook up a less trivial example where $\Gamma=G\times F$ where $G$ has a free action on a $\Lambda$--tree $Y$,
and $F$ is a free group acting freely on a simplicial tree $S$. $\Gamma$ acts on $S$ with kernel $G$.
Now take $Y_v=Y$ and $\phi_e=\Id$. Then $T_\calg=Y$, and $F$ is in the kernel of the action of $\Gamma$.
\end{rem*}

\begin{proof}
  If an element $g\in\Gamma$ fixes a point in $T_\calg$, then $g$ globally preserves the corresponding
equivalence class in $X$. Since this equivalence class has the structure of a tree with finite diameter, 
$g$ must fix a vertex in this equivalence class (there are no inversions because the action on $S$ has no inversion). 
Hence $g$ fixes a point of $X$, which means that $g$ fixes a point in a vertex tree.
\end{proof}

\subsection{Transverse coverings and graph of actions on $\bbR$--trees}
\label{sec_transverse_coverings}

In this section, we restrict to the case of a graph of actions on $\bbR$--trees along points.
We prove that an action on an $\bbR$--tree splits as such a graph of actions if
and only if it has a certain kind of covering by subtrees.
The argument could in fact be generalised to graph of actions on $\Lambda$--trees along points
but we won't need it.

\begin{dfn}[Transverse covering]\label{dfn_transverse_covering}\rm
Let $T$ be an $\bbR$--tree, and $(Y_u)_{u\in U}$ be a family of non-degenerate closed subtrees of $T$.
We say that $(Y_u)_{u\in U}$ is a transverse covering of $T$ if
\begin{description}
\item[\labelitemi\qua transverse intersection:] whenever $Y_u\cap Y_v$ contains more than one point, $Y_u=Y_v$;
\item[\labelitemi\qua finiteness condition:]  every arc of $T$ is covered by finitely many $Y_u$'s.
\end{description}
\end{dfn}

\begin{rem*}
  When $T$ is endowed with an action of a group $\Gamma$, we always require that the family $(Y_u)$ is $\Gamma$--invariant.
\end{rem*}

\begin{lem}\label{lem_transverse} Consider an action of a group $\Gamma$ on an $\bbR$--tree $T$.
If $T$ splits as a graph of actions on $\bbR$--trees along points $\calg$,
then the image in $T$ of non-degenerate vertex trees of $\calg$ gives a transverse covering of $T$.

Conversely, if $T$ has a $\Gamma$--invariant transverse covering, 
then there is a natural graph of actions $\calg$
whose non-degenerate vertex trees correspond to the subtrees of the transverse covering
and such that $T\simeq T_\calg$.
\end{lem}

\begin{proof}
We first check that the family of vertex trees of a graph of actions $\calg=(S,(Y_v),(\phi_e))$ forms
a transverse covering of $T_\calg$. We have already noted that vertex trees are closed in $T_\calg$.
The transverse intersection condition follows from the fact that edge trees are points.
To prove the finiteness condition, consider $x\in Y_u$ and $y\in Y_v$ and note that
$[x,y]$ is covered by the trees $Y_w$ for $w\in [u,v]$.

To prove the converse, we need to define the simplicial tree $S$ encoding the combinatorics of the gluing.
\begin{dfn}[Skeleton of a transverse covering] \label{dfn_skeletton}\rm
Consider a transverse covering $(Y_u)_{u\in U}$ of $T$.
The \emph{skeleton} of this transverse covering is the bipartite simplicial tree $S$ defined as follows:
\begin{itemize*}
  \item $V(S)=V_0(S)\cup V_1(S)$ where $V_1(S)=\{Y_u \,|\ u\in U\}$, and $V_0(S)$ is the set of points
of $T$ which belong to a least two distinct trees $Y_u\neq Y_v$
\item there is an edge between $x\in V_0(S)$ and $Y\in V_1(S)$ if and only if $x\in Y$.
\end{itemize*}
\end{dfn}

The connectedness of $S$ follows from the finiteness condition (using the fact that the subtrees $Y_u$ are closed in $T$).
Now let's prove the simple connectedness.
Consider a path $p=x_0,Y_0,x_1,\dots,x_{n-1},Y_{n-1},x_n$ in $S$, and let 
$\Tilde p=[x_0,x_1].[x_1,x_2]\dots [x_{n-1},x_n]$ be the corresponding path in $T$.
If $p$ does not backtrack, then $Y_i\cap Y_{i+1}=\{x_{i+1}\}$ so $\Tilde p$ does not backtrack.
Therefore, $x_0\neq x_n$ and $p$ is not a closed path.

Now there is a natural graph of actions $\calg$ corresponding to $S$:
for $x\in V_0(S)$, the corresponding vertex tree is the point $\{x\}$,
for $Y\in V_1(S)$, the corresponding vertex tree is $Y$,
and the gluing maps $\phi_e\co \{x\}\ra Y$ are given by inclusion.
Finally, consider the natural map $\Psi\co T_\calg\ra T$
given by the inclusion of vertex trees.
This application is an isometry in restriction to vertex trees,
and if $[a,b]$, $[a,c]$ are two arcs in $T_\calg$ lying in two 
distinct vertex trees $Y_1,Y_2$ of $T_\calg$, then $\Psi([a,b])\cap \Psi([a,c])\subset\Psi(Y_1)\cap \Psi(Y_2)$ 
is reduced to one point.
This implies that $\Psi$ is an isometry in restriction to each segment, and hence an isometry.
\end{proof}

\begin{rem*}
  We will often prefer using a transverse covering (or the graph of actions corresponding to such a covering)
to a general graph of actions because of the following \emph{acylindricity} property of the dual graph of actions
$\calg=(S,(Y_v)_{v\in V(S)},(\phi_e)_{e\in E(S)})$:
if two points $x\in Y_v$, and $x'\in Y_{v'}$ have the same image in $T_\calg$, then $v,v'$ are at distance
at most $2$ in $S$.
\end{rem*}

It is also worth noticing the following simple minimality result:
\begin{lem}
Consider an $\bbR$--tree $T$ endowed with a minimal action of $\Gamma$.
Consider $(Y_u)_{u\in U}$ an equivariant transverse covering of $T$, and let $S$ be
the skeleton of the transverse covering.

Then the action of $\Gamma$ on $S$ is minimal.
\end{lem}

\begin{proof}
  Assume that $S'\subset S$ is an invariant subtree.
Let $T'\subset T$ be the union of vertex trees of $S'$.
One easily checks that $T'$ is connected using the connectedness of $S'$.
Thus, by minimality of $T$, one has $T'=T$. 
Using the \emph{acylindricity} remark above, $S$ is contained in the $2$--neighbourhood of $S'$.
In particular, if $S'\neq S$, then $S$ contains a terminal vertex $v$.
By definition, every vertex in $V_0(S)$ has at least two neighbours, so $v\in V_1(S)$.
We thus get a contradiction since $Y_v$ is contained in $T'$, contradicting the transverse intersection property
of the transverse covering.
\end{proof}





\section{The action modulo infinitesimals, abelian d\'evissage}\label{sec_abelian}
In this section, we prove a weaker version of the cyclic d\'evissage theorem,
where (maybe non-finitely generated) abelian groups may appear in place of cyclic groups
(see Proposition \ref{abelian_devissage}).

Start with a finitely generated group $\Gamma$ acting freely on an $\bbR^n$--tree $T$ with $n\geq 2$,
and assume that $\Gamma$ is freely indecomposable.
Denote by $\bbR^{n-1}$ the maximal proper convex subgroup of $\bbR^n$,
and consider elements of $\bbR^{n-1}$ as \emph{infinitesimals}.
Now consider the $\bbR$--tree $\olT$ obtained from $T$ 
by identifying points at infinitesimal distance 
(this is often called the \emph{base change functor} in the literature, 
see for instance \cite{Chi_nontrivial}, \cite{JaZa_Chiswell},
\cite{Bass_non-archimedean}, see also \cite[Theorem 2.4.7]{Chi_book}).
Note that the canonical projection $f\co T\ra \olT$ preserves alignment, and that the preimage of a convex set
 is convex. The preimage in $T$ of a point of $\olT$ is thus an infinitesimal subtree of $T$.
Of course, the action of $\Gamma$ on $T$ induces an isometric action of $\Gamma$ on $\olT$.

However, this action generally fails to be free.
It may even happen that $\Gamma$ fixes a point $\ol{x}$ in $\olT$, but in this case, 
the d\'evissage theorem holds trivially 
since $\Gamma$ acts freely on the $\bbR^{n-1}$--tree $f\m(\ol{x})$.
Therefore, we assume that $\Gamma$ acts non-trivially on $\olT$, and, up to taking a subtree of $\olT$ and its preimage in $T$,
 we can assume that the action on $\olT$ is minimal, \ie that there is no non-empty proper invariant subtree.

 We first analyse how far from free this action can be.

\begin{fact} \label{fact_elem}
If a group $\Gamma$ acts freely on an $\bbR^{n}$--tree $T$, then the action of $\Gamma$ on $\olT$ satisfies the following: 
 \begin{itemize*}
  \item tripod fixators are trivial (a \emph{tripod} is the convex hull of 3 points which are not aligned)
  \item for every pair of commuting, elliptic elements $g,h\in \Gamma\setminus\{1\}$, 
$\Fix_{\olT} g=\Fix_{\olT} h$; in particular,  $\Fix_{\olT} g=\Fix_{\olT} g^k$ for $k\neq 0$;
  \item arc fixators are abelian; the global stabilizer of a line is maximal abelian if it is non-trivial;
  \item the action is superstable: for every non-degenerate arc $I\subset \olT$ with non-trivial fixator, for every
non degenerate sub-arc $J\subset I$, one has $\Stab I=\Stab J$.
  \end{itemize*}
\end{fact}

\begin{rem*}
  This fact does not use the fact that we have an $\bbR^n$--tree rather than a more general $\Lambda$--tree.
The statement holds for every $\Lambda$--tree $T$ with a free action of $\Gamma$ without inversion 
and every $\Lambda/\Lambda_0$--tree $\olT$ obtained from $T$
by killing a convex subgroup $\Lambda_0$ of infinitesimals.
\end{rem*}

\begin{proof}[Proof of the fact]
We start with the proof of the two first items.
  Consider $g\in\Gamma\setminus\{ 1\}$, and consider  $\Fix_\olT g$ its set of fix points in $\olT$.
The preimage of $\Fix_\olT g$ is the set of points in $T$ moved by an infinitesimal amount.
This set is either empty if $l_T(g)$ is not infinitesimal, or, it is the
set of points whose distance to the axis $A_g$ of $g$ is infinitesimal.
Therefore, if $\Fix_\olT g$ is not empty, 
then it is the image of $A_g$ in $\olT$ which contains no tripod since the quotient map preserves alignment. 
Moreover, for any element $h\in \Gamma$ commuting with $g$, $h$ globally preserves the axis of $g$, so $A_h=A_g$.
Therefore, if $h$ is elliptic in $\olT$ then it has the same set of fixed points as $g$.

Now we prove superstability and that arc fixators are abelian. 
Consider some non-degenerate arcs $\ol{J}\subset\ol{I}\subset \olT$, two elements
$g,h\in\Gamma\setminus \{1\}$ fixing pointwise $\ol{I}$ and $\ol{J}$ respectively. 
We want to prove that $h$ fixes $\ol{I}$ and commutes with $g$.
By hypothesis, the translation length of $g$ and $h$ are infinitesimal in $T$, 
their axes $A_g$ and $A_h$ must intersect in a subset of non-infinitesimal diameter. 

Now since the diameter of $A_g\cap A_h$ is much (infinitely) larger than $l_T(g)+l_T(h)$, $ghg\m h\m$
is elliptic in $T$ (see for instance \cite[Rk.~p.111]{Chi_book}). Since the action is free, this means that $g$ and $h$ commute
and in particular, the fixator of $\ol{I}$ is abelian.
This implies that $h(A_g)=A_g$, 
thus $A_h\supset A_g$ since $A_h$ is the maximal $h$--invariant linear subtree of $T$,
and $A_h=A_g$ by symmetry of the argument.
Therefore, $\Fix_\olT h=\Fix_\olT g$ and in particular, $h$ fixes $\ol{I}$.

Let's prove that the global stabilizer $\Gamma_l$ of a line $l\subset \olT$ is abelian.
Since $\Fix g=\Fix g^2$ for all $g\in \Gamma$, $\Gamma_l$ acts on $l$ by translations.
If the fixator $N_l$ of $l$ is trivial, then we are done. Otherwise, $N_l$ is a normal abelian subgroup of $\Gamma_l$,
and let $\Tilde l$ be its axis in $T$. Since $\Gamma_l$ normalizes $N_l$, $\Gamma_l$ preserves $\Tilde l$,
so $\Gamma_l$ acts freely by translations on $\Tilde l$, so $\Gamma_l$ is abelian.
Finally, any element normalising $\Gamma_l$ must preserve $l$, so $\Gamma_l$ is maximal abelian.
\end{proof}

Therefore, one can apply Sela's theorem which claims that superstable actions on $\bbR$--trees
are obtained by gluing equivariantly some simpler $\bbR$--trees along points (see definition \ref{dfn_goa}). 
In this statement a \emph{simplicial arc} in $\olT$
is an arc $[a,b]$ which contains no branch point of $\olT$ except maybe at $a$ or $b$.

\begin{thm}[Structure theorem {\cite[Theorem 3.1]{Sela_acylindrical}}, see also \cite{Gui_Rips}]\label{thm_structure}
  Let $(\olT,\Gamma)$ be a minimal action of a finitely generated group on
  an $\bbR$--tree.  Assume that $\Gamma$ is freely indecomposable,
that tripod fixators are trivial, and that the action is super-stable.
Then $\olT$ can be decomposed into a graph of actions on $\bbR$--trees along points, each vertex tree being either
\begin{enumerate*}
\item a point;
\item a simplicial arc, which is fixed pointwise by its global stabilizer;
\item a line $l$ together with an action $\Gamma_l\actson l$ having dense orbits, 
such that the image of $\Gamma_l$ in $\Isom(l)$ is finitely generated;
\item or an action on an $\bbR$--tree dual to an arational%
\footnote{A measured foliation on a surface with boundary is \emph{arational} if any non simply-connected leaf (or generalised leaf)
actually has a cyclic fundamental group and contains a boundary component of $\S$. Equivalently, $\calf$ is arational
if every simple closed curve having zero intersection with the measured foliation is boundary parallel.}
measured foliation on a $2$--orbifold (with boundary).
\end{enumerate*}
\end{thm}

\begin{rem*}
  In \cite{Gui_Rips}, simplicial arcs are incorporated in the skeleton of the decomposition of the action
(as edges of positive length) and hence do not appear in the statement of the theorem.

Since $\Gamma$ is torsion-free, the orbifold groups occurring in the structure theorem are actually surface groups.
\end{rem*}


\paragraph{Agglutination of simplicial arcs}
We now make the decomposition given in the structure theorem nicer
with respect to abelian groups.
In particular, we want to gather simplicial arcs having the same fixator into bigger
vertex subtrees. This will imply that the stabilizer of the new corresponding
vertex trees are maximal abelian. The goal is to reformulate the Structure Theorem as follows:

\begin{thm}[Reformulation of Structure Theorem]\label{thm_reformulation}
There is a $\Gamma$--invariant transverse covering of $\olT$ by a family $(\olY_u)_{u\in U}$ of non-degenerate closed subtrees
such that, denoting by $\Gamma_u$ be the global stabilizer of $\olY_u$, one of the following holds:
\begin{description}
\item[\labelitemi\qua abelian-type:] $\olY_u$ is an arc or a line, 
the image $\olG_u$ of $\Gamma_u$ in $\Isom(\olY_u)$ is finitely generated, and
$\Gamma_u$ is maximal abelian in $\Gamma$; moreover for any two  abelian-type subtrees $\olY_u\neq\olY_v$,
$\Gamma_u$ and $\Gamma_v$ don't commute;
\item[\labelitemi\qua surface-type:] or the action $\Gamma_u\actson \olY_u$ is dual to an arational measured foliation on a surface with boundary.
\end{description}
\end{thm}

\begin{proof}
Consider the transverse covering of $\olT$ by the non-degenerate vertex subtrees $(\olY_u)_{u\in U}$
of the decomposition given by the structure theorem \ref{thm_structure} (lemma \ref{lem_transverse}).
Consider the equivalence relation on $U$ generated by $u\sim u'$ if $\olY_u$ and $\olY_{u'}$
are simplicial arcs and the fixators of $\olY_u$ and $\olY_{u'}$ commute
(note that the fixators of these arcs are non-trivial since $\Gamma$ is freely indecomposable).
The commutation of the fixators implies that those fixators coincide 
since commuting elliptic elements have the same set of fix points (fact~\ref{fact_elem}).

For any equivalence class $[u]$, let $\olZ_{[u]}=\bigcup_{u\in[u]} \olY_u$.
We prove that $(\olZ_{[u]})_{u\in U/{\sim}}$ is the wanted transverse covering.
For $u$ such that $\olY_u$ is a simplicial arc, let $N_{[u]}$ be the common fixator of the simplicial arcs $\olY_u$ for $u\in [u]$. 
One has $\Fix N_u=\olZ_{[u]}$. As a matter of fact, $N_u$ cannot fix an arc in a surface type vertex tree since the fixator of an arc
in a surface-type vertex tree is trivial, and $N_u$ cannot fix an arc in a line-type vertex tree because tripod fixators are trivial.
This implies that $\olZ_{[u]}$ is closed and connected, and is a linear subtree of $\olT$ since tripod fixators are trivial.
Since the case of a semi-line is easy to rule out, $\olZ_{[u]}$ is either an arc of a line. 
In particular, the family of subtrees $(\olZ_{[u]})_{[u]\in U/\eqv} $ is a transverse covering of $\olT$.

The global stabilizer $\Gamma_{[u]}$ of $\olZ_{[u]}$ is maximal abelian in $\Gamma$:
this was already noted in the case where $\olZ_{[u]}$ is a line;
if $\olZ_{[u]}$ is an arc, then $\Gamma_{[u]}$ coincides with its fixator $N_u$
since there can be no reflection (because $\Fix g=\Fix g^2$) so in particular, $\Gamma_{[u]}$ is abelian.
Now any element $g$ commuting with the elements of
$N_u$ must globally preserve $\olZ_{[u]}$, so $g\in \Gamma_{[u]}$. Therefore, $\Gamma_{[u]}$ is maximal abelian.
If $u,v$ are such that $\Gamma_{[u]}$ and $\Gamma_{[v]}$ commute, 
then $N_u$ and $N_v$ commute, so $N_u=N_v$ and $[u]=[v]$.
\end{proof}

We now focus on the skeleton $S$ of this transverse covering, and we analyse the corresponding splitting of $\Gamma$. 
We prove that this splitting satisfies the abelian d\'evissage. We give a simple version before giving
a more detailed statement.

\begin{prop}[Abelian d\'evissage, simple version]\label{abelian_devissage_simple}
If a finitely generated freely indecomposable group $\Gamma$ acts freely on an $\bbR^n$--tree ($n\geq 2$), 
then $\Gamma$ can be written as the fundamental group of a finite graph of groups 
where 
\begin{itemize*}
  \item each edge group is abelian; more precisely, an edge group is either cyclic or fixes an arc in $\olT$;
  \item each vertex group acts freely on an $\bbR^{n-1}$--tree. 
\end{itemize*}
\end{prop}

\begin{prop}[Abelian d\'evissage, detailed version]\label{abelian_devissage}
If a finitely generated freely indecomposable group $\Gamma$ acts freely on an $\bbR^n$--tree, 
then $\Gamma$ can be written as the fundamental group of a finite graph of groups 
with 3 types of vertices named \emph{abelian}, \emph{surface} and \emph{infinitesimal}, and such that the following holds:
\begin{itemize*}
\item each edge is incident to exactly one infinitesimal vertex;
\item for each abelian vertex $v$, $\Gamma_v$ is abelian maximal in $\Gamma$, 
$\Gamma_v=\olG_v\oplus N_v$ where $\olG_v$ is a finitely generated (maybe trivial) free abelian group,
$N_v$ is an arc fixator, and the image in $\Gamma_v$ of all incident edges  coincide with the (maybe infinitely generated) abelian group $N_v$;
moreover, if $v\neq v'$ are distinct abelian vertices, then $\Gamma_v$ does not commute with any conjugate of $\Gamma_{v'}$;
\item for each surface vertex $v$, $\Gamma_v$ is the fundamental group of a surface $\S$ with boundary holding
an arational measured foliation; there is one edge for each boundary component of $\S$, and the image of its edge group in $\Gamma_v$
is conjugate to the fundamental group of the corresponding boundary component of $\S$; 
\item for each infinitesimal vertex $v$, $\Gamma_v$ acts freely on an $\bbR^{n-1}$--tree; moreover, any element $g\in\Gamma\setminus\{1\}$ commuting
with an element of $\Gamma_v\setminus\{1\}$ either belongs to $\Gamma_v$, or is conjugate into $\olG_w$
where $w$ is an abelian vertex neighbouring $v$.
\end{itemize*}
Finally, $\calg$ is $4$--acylindrical and any non-cyclic abelian subgroup of $\Gamma$ is conjugate into a vertex group.
\end{prop}

\begin{rem*}
A surface with empty boundary could occur in this graph of groups, but in this case, the graph
of groups contains no edge, and $\Gamma$ is a surface group.

Note that the edge and vertex groups could a priori be non-finitely generated in the abelian d\'evissage.
On the other hand, if one knew somehow\footnote{Note that the claim of Remeslennikov that limit groups act freely on a $\Lambda$--tree
with $\Lambda$ finitely generated would imply that abelian subgroup are finitely generated since they are isomorphic 
to subgroups of $\Lambda$ since an abelian subgroup of $\Gamma$ acts freely by translation on its axis, and is 
thus isomorphic to a subgroup of $\Lambda$.}
that abelian subgroups of $\Gamma$ were finitely generated, then the finite presentation of
$\Gamma$ would follow easily: finite generation of edge groups would imply the finite generation of vertex groups (since $\Gamma$ is finitely generated),
and thus vertex groups would be finitely presented by induction hypothesis.
\end{rem*}

If one knows that arc fixators of $\olT$ are cyclic, then it is immediate to deduce the conclusion of the cyclic d\'evissage theorem 
from the abelian d\'evissage theorem. The strategy for the proof of the Cyclic D\'evissage Theorem will thus consist
in finding an $\bbR^n$--tree $T'$ such that arc fixators of $\olT'$ are cyclic (see next section).

\begin{proof}[Proof of the simple version from the detailed version]
  The claim about edge groups follows from the fact that each edge is either incident on a surface vertex or
on an abelian vertex. The claim about vertex groups follows from the fact that countable torsion free abelian groups and surface groups
holding an arational measured foliation (which are free groups in the case where the surface have non-empty boundary)
have a free action on an $\bbR$--tree.
\end{proof}

\begin{proof}[Proof of the detailed version]
Let $S$ be the skeleton of the transverse covering given by the reformulation of the structure theorem (Theorem  \ref{thm_reformulation}).
We prove that the graph of groups decomposition $\Gamma=\pi_1(\calg)$ induced by the action of $\Gamma$ 
on $S$ satisfies the abelian d\'evissage theorem.

Remember that $S$ is bipartite, with $V(S)=V_0(S)\dunion V_1(S)$ where $V_1(S)$ is the set of non-degenerate subtrees
in the transverse covering, and $V_0(S)$ is the set of points of $\olT$ which belong to at least two distinct subtrees
of the transverse covering. The set $V_0(S)$ will be the set of our infinitesimal vertices.
Since the stabilizer of such a vertex fixes a point in $T$, it acts freely on an $\bbR^{n-1}$--tree.
By the reformulation of the Structure Theorem, $V_1(S)$ is a disjoint union of abelian-type and surface type vertices
$V_1(S)=V_{ab}\dunion V_{surf}$,
where $V_{ab}$ is the set of vertices corresponding to abelian-type subtrees and $V_{surf}$ to surface-type subtrees (excluding tori).
The fact that $S$ is bipartite means that each edge of $\calg$ is incident on exactly one infinitesimal vertex.

Let's first consider an abelian vertex $v$, and let $N_v$ be the fixator of the linear subtree $\olY_{v}$, and $\olG_v$
the image of $N_v$ in $\Isom(\olY_v)$. The direct sum follows from the fact that 
$\Gamma_v$ is abelian and that $\olG_v$ is a free abelian group.
The only thing to check is that the image of all edge groups incident on $v$
 coincide with $N_v$. This follows from the fact that an edge $e\in E(S)$
is a pair $(x,\olY_v)$ where $x\in\olY_v$, so the fixator of $e$ is the stabilizer of $x$ in $\Gamma_v$,
which is $N_v$.

The acylindricity follows from the fact that if two edges $(x,\olY_v)$ $(x',\olY_{v'})$ have 
commuting fixators $\Gamma_e,\Gamma_{e'}$, then $\Gamma_e$ and $\Gamma_{e'}$ have the same (non-empty) set of
fix points in $\olT$ which is either a point, or an abelian subtree. In the first case, one has $x=x'$
so the two edges have a common endpoint at an infinitesimal vertex. In the second case,
the endpoints $x$ and $x'$ of the two edges are at distance at most $2$ in $S$.
The acylindricity implies that any non-cyclic abelian group $A$ is conjugate into a vertex group
since otherwise, a non-trivial subgroup of $A$ would fix its axis in $S$, contradicting acylindricity.

Let's turn to a surface vertex $v$. We know that its stabilizer $\Gamma_v$ 
is the fundamental group of a surface with boundary holding an arational measured foliation.
Moreover, since $\Gamma$ is freely indecomposable,
edge fixators are non-trivial. The fixator of an edge $e=(x,\olY_v)$
is non-trivial and fixes a point in $\olY_v$. Since $\olY_v$ is dual to an arational measured
foliation on a surface $\Sigma$,  
the elliptic elements of $\Gamma_v$ are exactly those which can be conjugate into the fundamental group of a boundary component of $\Sigma$,
and an elliptic element fixes exactly one point.
Thus, $\Gamma_e$ is conjugate to a boundary component $b_e$ of $\Sigma$ (and not to a proper subgroup since $\Gamma_e$ is the whole
stabilizer of $x$ in $\Gamma_v$).
Moreover, if two edges $e=(x,\olY_v),e'=(x',\olY_v)$ of $S$ correspond to the same boundary component $b_e=B_{e'}$ of $\Sigma$, then
$\Gamma_e$ and $g\Gamma_{e'}g\m$ for some $g\in \Gamma_v$, which implies that $x'=g.x$ since $\Gamma_e$ fixes
exactly one point in $\olY_v$, so $e'=g.e$. 
This proves that two distinct edges of $\calg$ incident on a surface vertex of $\calg$
correspond to \emph{distinct} boundary components of the surface.
If a boundary component of $\Sigma$ does not correspond to any incident edge,
it is easy to check that $\Gamma$ has a non-trivial free splitting, contradicting the hypothesis.

There remains to check the last affirmation about elements commuting with an element stabilising
an infinitesimal vertex. So let $v\in V_0(S)$,  $g\in \Gamma_v\setminus\{1\}$ and $h\in \Gamma$
commuting with $g$. 
If $h$ is elliptic, then $h\in \Gamma_v$ since two commuting  elliptic elements have the same fixed points.
If $h$ is hyperbolic, then $g$ must fix pointwise its axis $A_h$. Since surface-type subtrees have trivial arc fixators,
 $A_h$ cannot meet any surface-type subtree in more than a point.
So $A_h$ is contained in a union of abelian subtrees. But since $g$ fixes $A_h$,
$A_h$ is a single abelian subtree. Now $x\in A_h$ since $g$ fixes no tripod and the last claim follows.
\end{proof}



\section{Obtaining cyclic arc fixators}\label{sec_flawless}

To prove the cyclic d\'evissage theorem, we will find an $\bbR^n$--tree $T'$
such that arc fixators of the $\bbR$--tree $\olT'$ are cyclic. The cyclic d\'evissage theorem
will then follow from the abelian d\'evissage theorem.

\begin{thm}\label{thm_cyclic_flawless}
  Assume that a freely indecomposable finitely generated group $\Gamma$ acts freely
on an $\bbR^n$--tree $T$. Then $\Gamma$ has a free action on an $\bbR^n$--tree $T'$
such that the action on the $\bbR$--tree $\olT'$ obtained by killing infinitesimals
has cyclic arc fixators.
\end{thm}

The strategy of the proof is the following: starting with the 
decomposition of $\Gamma$ given by the structure theorem for the action of $\Gamma$ on $\ol T$,
we build a more adapted decomposition of $\Gamma$ as a graph of groups with abelian edge groups, 
such that the action of its non-abelian vertex groups
on their minimal subtree in $\overline{T}$ has cyclic arc fixators. Taking the Bass--Serre tree of this decomposition as a skeleton,
we define a graph of actions on $\bbR^n$--trees satisfying the conclusion of Theorem \ref{thm_cyclic_flawless}.

But to clarify the exposition, we first sketch the proof of Theorem \ref{thm_cyclic_flawless} 
in two important particular cases.

\subsection{Main examples}

\subsubsection{Acylindrical case}

The first case concerns a decomposition of $\Gamma$ as an amalgamated product
of non-abelian groups over a maximal abelian subgroup. The CSA property of $\Gamma$ 
immediately implies the $1$--acylindricity of the splitting.

\begin{lem}
Consider a free action of a finitely generated group $\Gamma$ on an $\bbR^n$--tree, and assume that
$\olT$ is a simplicial tree, dual to an amalgamated product $\Gamma=A*_C B$.
In particular, $A$ and $B$ both stabilize a $\bbR^{n-1}$--tree in $T$.
For simplicity, also assume that $C$ contains an element whose translation length has magnitude $n-1$. 

If $C$ is not cyclic and maximal abelian in both $A$ and $B$, 
then $\Gamma$ has a free action on an $\bbR^{n-1}$--tree $T'$.  
\end{lem}

\begin{rem*}
  The hypothesis that $C$ is not cyclic is crucial here.
Indeed, a group of the from $\Gamma=F_2*_C F_2$ where $C$ is maximal cyclic
in both copies of $F_2$ has a free action on a $\bbZ^2$--tree (see \cite{Bass_non-archimedean}).
However, if $C$ is chosen so that $\Gamma$ is not a free product of surface or abelian groups,
then $C$ has no free action on any $\bbR$--tree.
\end{rem*}

\begin{proof}
Take $S=\olT$ endowed with its natural set of vertices and edges.
For each vertex $v\in V(S)$, let $Y_v$ be the preimage of $v$ in 
the $\bbR^n$--tree $T$. Note that $Y_v$ is an $\bbR^{n-1}$--tree.
For each edge $e\in E(S)$, denote by $\Gamma_e$ its fixator, 
and let $\lambda_e$ be the axis of $\Gamma_e$ in $Y_{t(e)}$.
We will construct $T'$ as a graph of actions on $\bbR^{n-1}$--trees 
with $S$ as skeleton, and whose vertex and edge trees are $Y_v$ and $\lambda_e$.
The main point is to define the gluing isometries
so that the resulting action is free.

One technical problem which may occur in general is that $\lambda_e$ may not be isometric to $\lambda_{\ol e}$\,: 
it might be isometric
to a strict convex subset of $\bbR^{n-1}$. However, our hypothesis which claims that some element of 
$C$ has a translation length of magnitude $n-1$ prevents this.
In general, this problem will be fixed by enlarging the vertex trees using an \emph{end completion} procedure.

Now, we need to define the gluing maps $\phi_e$. To prove that we get a free action, we will apply 
our criterion for freeness (lemma \ref{lem_criterion_free}).
Indeed, we will choose the maps $\phi_e$ so that no point of $Y_v$ gets identified with a point of $Y_w$ where
$v$ and $w$ are at distance more than $2$ in $S$.
This is of course automatic if for each vertex $v\in V(S)$ and for each pair of distinct 
edges $e,e'$ incident on $v$, $\lambda_e$ does not intersect $\lambda_{e'}$.

We argue that for two edges $e\neq e'$ incident on $v$, $\lambda_e\cap\lambda_{e'}$
 is small: it has magnitude at most $n-2$. This is where the non-cyclicness of $C$ is crucial.
The hypothesis that $C$ is maximal abelian in both $A$ and $B$ implies that 
for each pair of distinct edges, $\Gamma_e$ and $\Gamma_{e'}$ don't commute
 (this is an easy exercise about the CSA property).
The smallness of  $\lambda_e\cap\lambda_{e'}$ 
follows from the following fact applied to the action of $\Gamma_v$ on $Y_v$ with $p=q=n-1$, $H=\Gamma_e$, and $H'=\Gamma_{e'}$. 

\begin{fact}[Infinitesimal intersection of axes of non-cyclic groups]\label{fact_small}
  Let $\Gamma\actson Y$ be a free action of a group on an $\bbR^p$--tree.
Let $H,H'$ be two non-commuting abelian subgroups of $\Gamma$. 
Assume that for some $q\leq p$, the subgroups of $H$ and $H'$
consisting of elements whose translation length is of magnitude at most $q$
are both non-cyclic.

Then the intersection of the axes of $H$ and $H'$ has magnitude at most $q-1$.
\end{fact}


We keep the proof of the fact for later.
\begin{dfn}\label{dfn_branching_locus}\rm
The \emph{branching locus of axes} in $Y_v$ is the set 
$$B_v=\bigcup_{\text{\textup{$e\neq e'\in E(S)$ incident on $v$}}} \big[\lambda_e\cap \lambda_{e'}\big].$$
\end{dfn}

The fact says that this branching locus is a countable union of sets of magnitude at most $n-2$.

Now choose an oriented edge $e$, and any isometry $\phi_e\co \lambda_{\ol e}\ra \lambda_e$.
Up to composing $\phi_e$ by a generic translation of $\lambda_e\simeq\bbR^{n-1}$,
we can assume that $\phi_e(B_{o(e)})$ does not intersect $B_{t(e)}$ (there are countably many
classes of prohibited translations mod $\bbR^{n-2}$).
Now extend this choice of $\phi_e$ equivariantly to get an equivariant graph of actions on $\bbR^{n-1}$--trees.

There remains to prove that no point of $Y_v$ gets identified with a point of $Y_w$ where
$v$ and $w$ are at distance more than $2$ in $S$.
Take $x\in Y_v$ and assume that $x$ is identified with two points $y\in Y_u$ and $z\in Y_w$
where $u,w$ are two distinct neighbours of $v$. If $e,e'$ denote the edges joining $v$ to $u$ and $w$,
we have that $x\in\lambda_e\cap\lambda_{e'}$ so $x\in B_v$. 
Therefore, our choice of gluing isometries implies that $y\notin B_u$ and $z\notin B_w$. This implies that $y$ and $z$ are not identified
with any other point. Therefore, Lemma \ref{lem_criterion_free} applies and the action dual to our graph of actions is free.
\end{proof}

\begin{proof}[Proof of Fact \ref{fact_small}]
  This is just a consequence of the fact that if the axes of two hyperbolic elements $g,h$ 
intersect in a segment whose diameter is larger than the sum of the translation lengths of $g$ and $h$, 
then the commutator $[g,h]$
fixes a point, so that $g$ and $h$ commute.

Up to taking subgroups, one can assume that every element of $H$ and $H'$
have translation length of magnitude at most $q$.
We just need to prove that for any positive $\eps\in\bbR^q\setminus\bbR^{q-1}$, 
$H$ (resp.\ $H'$) contain non-trivial elements of translation length at most $\eps$. 
This clearly holds if some element of $H\setminus\{1\}$ has a translation length of magnitude at most $q-1$.
Otherwise, consider a morphism $\rho\co H\ra\bbR^q$ having the same translation length function.
Composing $\rho$ by the collapse of infinitesimals in $\bbR^q$, we get an embedding of $H$ as a subgroup of $\bbR$.
This subgroup has to be dense since $H$ is not cyclic.
\end{proof}

\subsubsection{Non acylindrical case}

The next case occurs when one allows different edge fixators to commute.
In this case we need to modify the graph of groups to get a more adapted one.
The major example is the following: $\Gamma=A*_C$ where both inclusions of $C$ in $A$ are identical.
Then we can consider the following modified decomposition: 
$\Gamma=A*_C (C\oplus \bbZ)$.
The corresponding Bass--Serre tree has the following nice feature which will be generalised later:
if $\Gamma_e$ commutes with $\Gamma_{e'}$ then $e,e'$ are incident on a common vertex whose
stabilizer is abelian (see Lemma \ref{lem_adapted_decompo}).

\begin{lem}
Consider a free action of a group $\Gamma$ on an $\bbR^n$--tree $T$.
Assume that $\Gamma$ splits as an amalgamated product $A*_C \Hat C$ where
$A$ is non-abelian and preserves an $\bbR^{n-1}$--subtree of $T$, 
that $\Hat C$ is abelian and $\Hat C\simeq C\oplus \ol C$ 
for some non-trivial finitely generated free abelian group $\ol C$.
Also assume for simplicity that $C$ contains an element whose translation length has magnitude $n-1$.

If $C$ is not cyclic, then
$\Gamma$ has a free action on an $\bbR^{n-1}$--tree.
\end{lem}

\begin{proof}
We want to define a $\Gamma$--equivariant graph of actions on $\bbR^{n-1}$--trees.
Let $S$ be the Bass--Serre tree of the splitting $\Gamma=A*_C \Hat C$.
We denote by $a$ and $c$ the vertices of $S$ whose stabilizers are $A$ and $\Hat C$ respectively,
and by $e_0$ the edge connecting them with $t(e_0)=c$.

For each $v\in V(S)$ in the orbit of $a$, define $Y_v=\Tmin(\Gamma_v)$
as the minimal $\Gamma_v$--invariant subtree of $T$.
Remember that $B_a$ denote the branching locus of axes of $Y_a$, \ie the set
of points of $Y_a$ which belong to the intersection of the axes of two distinct conjugates of $C$.
Let $D=\{d(x,x')|\ x,x'\in B_a\}$  be the set of mutual distances between points
of the branching locus of axes. By Fact \ref{fact_small}, 
$D$ is a countable union of sets of magnitude at most $n-2$ (as a metric space);
in other words, each of these countably many sets has the property that
the difference of two of is elements lies in $\bbR^{n-2}$.

For each $v\in \Gamma.c$, define $Y_v$ as a copy of $\bbR^{n-1}$.
We choose any action of $\Hat C=C\oplus \ol C$ on $Y_v$ so that the translation length
of elements of $C$ is the same as in $T$, and so that the translation length
of any element of $\Hat C\setminus C$ lies in $\bbR^{n-1}\setminus D$.
We extend equivariantly this action to an action of $\Gamma$ on the union of those copies of $\bbR^{n-1}$.

Define $\lambda_{e_0}$ as $Y_c=\bbR^{n-1}$ and $\lambda_{\ol e_0}$ as the axis of $C$ in $Y_a$.
Since $C$ contains an element of translation length $n-1$, $\lambda_{\ol e_0}\simeq \bbR^{n-1}$,
and choose for $\phi_{e_0}$ any isometry.
Define $\lambda_e$ and $\phi_e$ for every edge by equivariance.

Thus, we have defined a graph of actions on $\bbR^{n-1}$--trees, and we need to prove
that the action of $\Gamma$ on the dual $\bbR^{n-1}$--tree is free.
To apply our criterion for freeness (Lemma \ref{lem_criterion_free}), 
there remains to prove that no point of $Y_v$ gets identified with a point of $Y_w$ where
$v$ and $w$ are at distance more than $4$ in $S$.

Assume on the contrary that there are two such points at distance $5$. 
Denote by $u=u_0,u_1,\dots,u_5=v$ the vertices of $[u,v]$, $e_i$ the edge $[u_{i-1},u_i]$
and let $x_i\in Y_{u_i}$ such that $\phi_{e_i}(x_{i-1})=x_i$.
If $u_i\in\Gamma.a$ for some $i\in\{1,\dots,4\}$, then $x_i$ lies in $B_{v_i}$.
Since any edge of $S$ joins a vertex in $\Gamma.a$ to a vertex in $\Gamma.c$,
one can find two indices $i,i+2\in\{1,\dots,4\}$ such that $u_i,u_{i+2}\in \Gamma.a$.
By equivariance, we can assume that $u_i=a$, $u_{i+1}=c$, and $u_{i+2}=g.a$
for some $g\in \Hat C$.
Since $u_{i+2}\neq u_i$, $g\notin C$.
The two points $x_i$ and $g\m x_{i+2}$ both lie in $B_a$,
so $d_{Y_a}(x_i,g\m x_{i+2})\in D$. 
Besides, $\phi_{\ol e_{i+2}}(x_{i+2})=x_{i+1}$ 
so $\phi_{e_{i+1}}(g\m x_{i+2})=g\m x_{i+1}$.
Since $\phi_{e_{i+1}}(x_i)=x_{i+1}$,
one gets $d_{Y_c}(x_{i+1},g\m x_{i+1})=d_{Y_a}(x_i,g\m x_{i+2}) \in D$,
a contradiction with our choice for the action of $\Hat C$ on $Y_c$
since $g\notin C$.
\end{proof}

\subsection{General case}

The following definition will be convenient.
\begin{dfn}[Flawless $\bbR^n$--trees]\rm
  Let $T$ be an $\bbR^n$--tree endowed with an action of $\Gamma$, and let $\olT$ be the 
$\bbR$--tree  obtained by killing infinitesimals.
One says that $T$ is \emph{flawless} if $\olT$ has cyclic arc fixators.
\end{dfn}

By extension, we will also say that the $\bbR$--tree $\olT$ is flawless accordingly.

With this terminology, Theorem \ref{thm_cyclic_flawless} means that 
any finitely generated group having a free action on an $\bbR^n$--tree,
has a \emph{flawless} free action on another $\bbR^n$--tree.

\subsection{Gluing flawless trees along infinitesimals is flawless}

We will build our flawless action as a graph of flawless actions over infinitesimal subtrees.
We first check that such an action is always flawless.

\begin{lem}\label{lem_flawless1}
Consider an action of a group $\Gamma$ on an $\bbR$--tree $\ol T$.
Assume that $\ol T$ is transversely covered (Definition \ref{dfn_transverse_covering}) 
by a $\Gamma$--invariant family of subtrees $(\olY_v)$
such that the action on $\olY_v$ of its global stabilizer is flawless.

Then $\Gamma\actson \ol T$ is flawless.
\end{lem}

\begin{proof}
Consider an arc $\ol I\subset \olT$. Using the finiteness condition in a transverse covering, 
up to making $\ol I$ smaller, one can assume that $\ol I$ is contained in a
non-degenerate subtree $\olY_v$ of the transverse covering. Any element fixing $\ol I$ must preserve $\olY_v$ by the transverse intersection
property. Since $\Gamma_v\actson\olY_v$ is flawless, this implies that the fixator of $\ol I$ is cyclic, which means that $\ol T$ is flawless.  
\end{proof}

We immediately deduce the following lemma:
\begin{lem}[Gluing flawless trees is flawless]\label{lem_flawless2}
  Consider a graph of flawless actions on $\bbR^n$--trees over infinitesimal subtrees.

Then the $\bbR^n$--tree dual to this graph of actions is flawless.
\end{lem}

Recall that a graph of actions on $\bbR^n$--trees is \emph{over infinitesimal subtrees}
if the edge subtrees are infinitesimal.

\begin{proof}
Let $S$ be the skeleton of the graph of actions, and denote by $Y_v$ the vertex $\bbR^n$--trees
and by $\phi_e\co \lambda_{\ol e}\ra\lambda_e$
the gluing isometries between infinitesimal closed subtrees of the corresponding vertex trees.
Let $T$ be the $\bbR^n$ tree dual to this graph of actions, and let $\olT$ be the $\bbR$--tree obtained by killing infinitesimals.
  The images $\olY_v$ of the vertex trees in $\olT$ give a transverse covering of $\olT$ 
because the gluing occurs along infinitesimal trees. 
\end{proof}

\subsubsection{Decomposition of the group}
The next step in the proof is to build
 an adapted decomposition of $\Gamma$  from the structure of the
action of $\Gamma$ on $\ol T$. This is the analogue of the transformation of $\Gamma=A*_C$ into $\Gamma=A*_C (C\oplus\bbZ)$
in the second example.
We state it in terms of the corresponding Bass--Serre tree.
 
\newcommand{\Vab}{V_{\mathrm{ab}}}
\renewcommand{\theenumi}{\roman{enumi}}
\renewcommand{\labelenumi}{\rm(\theenumi)}

\begin{lem}\label{lem_adapted_decompo}
  Let $\Gamma$ be a finitely generated group acting freely on an $\bbR^n$--tree $T$.
Then $\Gamma$ has an action on a simplicial tree $S$ such that
$S$ is bipartite for the partition of $V(S)$ into the set of vertices $\Vab(S)$ whose stabilizer is abelian
and the set $V_F(S)$ (F for flawless) whose stabilizer is non-abelian and such that
\begin{enumerate}
\item\label{enum_nonabelian} $\forall v\in V_F(S)$, the action of $\Gamma_v$ on the minimal subtree $\Tmin(\Gamma_v)$
is flawless;
\item\label{enum_abelian} $\forall v\in \Vab(S)$, $\Gamma_v$ is maximal abelian in $\Gamma$, 
and $\Gamma_v=N_v\oplus \ol\Gamma_v$ where $N_v$ coincides with the fixator of all edges incident on $v$
and $\olG_v$ is a finitely generated free abelian group;
\item\label{enum_edge} $\forall e\in E(S)$, $\Gamma_e$ is abelian, not cyclic, and the translation length of its elements in $T$
is infinitesimal;
\item\label{enum_commutation} if $e,e'\in E(S)$ are distinct edges such that $\Gamma_e$ and $\Gamma_{e'}$ commute, then $e$ and $e'$ share a common vertex  in $\Vab(S)$.
\end{enumerate}
\end{lem}

\begin{proof}
  Start with $S_0$ the tree dual to the graph of groups given by the Abelian D\'evissage Theorem \ref{abelian_devissage}.
We can assume that the action $\Gamma\actson S_0$ is minimal.
We keep the terminology of Proposition \ref{abelian_devissage} about abelian type vertices:
by an \emph{abelian type vertex} we mean a vertex of $S_0$ of abelian type in the sense of Proposition \ref{abelian_devissage};
we will not use this term to designate a vertex of $S_0$ (or $S$) whose stabilizer is abelian (\ie an element of $\Vab(S)$).

Take $S$ the tree obtained from $S_0$ by collapsing all the edges whose fixator is cyclic.
Denote by $\pi\co S_0\ra S$ the collapsing map.
Most requirements for $S$ will follow from the properties of $S_0$ described in Proposition \ref{abelian_devissage}.

First, for each edge $e\in E(S_0)$, the translation length of $\Gamma_e$ is infinitesimal because 
because one of the vertices of $e$ is of infinitesimal type. Therefore, the same holds for
edge fixators of $S$ and point (\ref{enum_edge}) is clear.

Note that all edges of $S_0$ incident on a surface type vertex are collapsed by $\pi$,
so any edge of $S_0$ which is not collapsed under $\pi$ joins an infinitesimal type vertex $u_0$ to an abelian type vertex $v_0$
whose stabilizer is of the form $\Gamma_{v_0}=N_{v_0}\oplus\ol\Gamma_{v_0}$ for some non-cyclic $N_{v_0}$.
Since all the edges incident on $v_0$ have the same fixator $N_{v_0}$, no edge incident on $v_0$ is collapsed.
Denoting $v=\pi(v_0)$, it follows that $\pi\m(v)=\{v_0\}$ so $\Gamma_v=\Gamma_{v_0}$ and $v$ satisfies
(\ref{enum_abelian}). 
To complete the proof of (\ref{enum_abelian}), we only need to check that for all $v$ with $\Gamma_v$
abelian, $\pi\m(v)=\{v_0\}$ for an abelian-type vertex as above.
Let $v\in V(S)$ with $\Gamma_v$ abelian. 
The subtree $\pi\m(v)$ cannot contain any surface-type vertex (because its stabilizer is non-abelian).
It can neither contain a vertex neighbouring a surface-type vertex (because the corresponding edge would be collapsed in $S$).
If $\pi\m(v)$ contains an infinitesimal-type vertex $v_0$, then $\Gamma_{v_0}$ is abelian, and all its neighbours
are of abelian type. There are at least two such neighbours $u_1,u_2$ by minimality, and $\Gamma_{u_1}$ and $\Gamma_{u_2}$ 
commute because $\Gamma_{v_0}$ is abelian (commutative transitivity).
Since $\Gamma_{u_1}$ is maximal abelian in $\Gamma$, CSA property implies that $u_1$ and $u_2$ cannot be in the same orbit.
But Proposition \ref{abelian_devissage} prohibits the commutation of stabilizers of abelian-type vertices 
in distinct orbits, a contradiction.
It follows that the subtree $\pi\m(v)$ is reduced to a single abelian type vertex, and 
(\ref{enum_abelian}) is proved.

It follows that $S$ is bipartite: consider $e\in E(S)$ and $\Tilde e\in E(S_0)$ its preimage;
since $\Gamma_{\Tilde e}$ is non-cyclic, $\Tilde e$ has no surface type endpoint, 
so $\Tilde e=[x,a]$ for an infinitesimal vertex $x$ and an abelian type vertex $a$;
since all edges incident on $a$ have the same non-cyclic fixator, $\pi\m(\pi(a))=a$ so
the stabilizer of $\pi(a)$ is abelian; the stabilizer of $\pi(x)$ can't be abelian by the
argument above.

Now take $v\in V(S)$ with $\Gamma_v$ non-abelian, and let's prove that the action of $\Gamma_v$ on $\Tmin(\Gamma_v)$
is flawless. Remember that $S_0$ was obtained as the skeleton of the transverse covering of $\ol T$ given by 
Theorem \ref{thm_reformulation}, so to each vertex $v_0\in V(S_0)$ corresponds either a point or a non-degenerate subtree
$\olY_{v_0}$ (according to whether $v_0$ is of infinitesimal type or not). 
Denote by $Y_{v_0}$ the preimage of the $\bbR$--tree $\olY_{v_0}$ in the $\bbR^n$--tree $T$.
Note that the action of $\Gamma_{v_0}$ on $Y_{v_0}$ is flawless except when $v_0$ is of abelian type with $N_{v_0}$ non-cyclic,
which occurs only when $\Gamma_{\pi(v_0)}$ is abelian.
Now consider $Y_v=\bigcup_{v_0\in \pi\m(v)} Y_{v_0}$, and let $\olY_v$ its image in $\ol T$. 
If $\olY_v$ is reduced to a point, then $\Gamma\actson Y_v$ and thus $\Gamma_v\actson\Tmin(\Gamma_v)$
are clearly flawless.
Otherwise, $\olY_v$ it is transversely covered by the non-degenerate subtrees of the family $(Y_{v_0})_{v_0\in \pi\m(v)}$. 
Moreover, each action $\Gamma_{v_0}\actson Y_{v_0}$ is flawless, so the action of $\Gamma_v$ on $Y_v$ is flawless by Lemma \ref{lem_flawless1}.
Since the flawless condition is stable under taking invariant subtrees, this proves (\ref{enum_nonabelian}).

Now consider two edges $e,e'\in E(S)$ whose stabilizers commute.
Denote by $\Tilde e,\Tilde e'\in E(S_0)$ their preimages.
Since $\Gamma_{\Tilde e}$ is non-cyclic, its non-infinitesimal endpoint is an abelian type vertex $u$.
Similarly, $\Tilde e'$ has an abelian type endpoint $u'$.
By commutative transitivity, $\Gamma_u$ commutes with $\Gamma_{u'}$.
Proposition \ref{abelian_devissage} then implies that $u=u'$, so $\Tilde e$ and $\Tilde e'$
share an abelian type vertex. Since $\pi(u)$ has an abelian stabilizer, (\ref{enum_commutation})
is proved.

This concludes the proof of the lemma.
\end{proof}

\subsubsection{Construction of the graph of actions}
Now we use the tree $S$ constructed in the previous lemma as a skeleton for a graph of actions on $\bbR^n$--trees
to build a flawless action of $\Gamma$ on an $\bbR^n$--tree.

Let's first define the vertex $\bbR^n$--trees $(Y_v)_{v\in V(S)}$ of our graph of actions.
For $v\in \Vab(S)$, we take for $Y_v$ a copy of $\bbR^{n-1}$
(we will define the action of $\Gamma_v$ on $Y_v$ later).
For $v\in V_F(S)$, we would like to take $Y_v=\Tmin(\Gamma_v)$, the minimal $\Gamma_v$--invariant
subtree of $T$.
For each edge $e\in E(S)$, denote by $l_e$ the axis of $\Gamma_e$ in $\Tmin(\Gamma_{o(e)})$
(resp.\ $l_e=\bbR^{n-1}$ if $t(e)\in \Vab(S)$).
We would next like to glue $l_e$ on $l_{\ol e}$ for each edge $e$.
Unfortunately, this might not be possible as it might happen that $l_e$ is not isometric to $l_{\ol e}$.

However, for every vertex $v\in V_F(S)$, and for any edge $e$ incident on $v$,
$l_e$ has magnitude at most $n-1$ because $\Gamma_e$ is non-cyclic and $\Tmin(\Gamma_{o(e)})$ is flawless.
As a remedy, we are going enlarge the vertex trees so that all axes of edge groups become
isometric to $\bbR^{n-1}$.

In what follows, we call \emph{line} in an $\bbR^n$--tree $T$ a \emph{maximal} linear subtree of $T$.

\begin{fact}[End completion, {\cite[Appendix E]{Bass_non-archimedean}}]
If $(Z,\Gamma)$ is an action on $\bbR^n$--tree, and $(l_e)$ is an invariant family of 
lines in $Z$ with magnitude at most $n-1$, 
then there is a natural enlargement $\Hat Z$ of $Z$ (endowed with an action of $\Gamma$) 
such that each $l_e$ is contained in a unique maximal line $\Hat l_e$ of $\Hat Z$,
and $\Hat l_e$ is isometric to $\bbR^{n-1}$. 
\end{fact}

\begin{proof}
  This fact follows for instance from \cite[Appendix E]{Bass_non-archimedean}:
take $\Hat Z$ to be the $\bbR^{n-1}$--neighbourhood of $Z$ in the $\bbR^{n}$--fulfilment of $Z$.
We  give an alternative simple sketch of proof for completeness, under the assumption 
that any two distinct lines of the family intersect in a segment (and not a semi-line for example).
This assumption is satisfied in our setting.

Fix a line $l_e$, choose an embedding $j_e$ of $l_e$ into $\bbR^{n-1}$, and
glue $Z$ to $\bbR^{n-1}$ along the maximal line $l_e$ using $j_e$. There is actually no choice for doing this
since any two embeddings of $l_e$ into $\bbR^{n-1}$ differ by an isometry of $\bbR^{n-1}$.
It is easily seen that the glued tree $Z\cup_{j_e}\bbR^{n-1}$ is an $\bbR^n$--tree
(although $j_e(l_e)$ is generally not  closed in $\bbR^{n-1}$ in the sense of $\Lambda$--trees).
The additional assumption we made says
that any other line $l_{e'}$ is a maximal linear subtree in the extended tree.
Therefore, we can iterate this construction, 
and the obtained tree does not depend on the order chosen to extend the lines.
Since an increasing union of $\bbR^n$--trees is an $\bbR^n$--tree, the fact is proven.
\end{proof}

This allows us to define our vertex and edge trees as follow:
\begin{itemize*}
\item for $v\in \Vab(S)$, take for $Y_v$ a copy of $\bbR^{n-1}$;
\item for $v\in V_F(S)$, take $Y_v=\widehat{\Tmin(\Gamma_v)}$ 
to be the end completion of $\Tmin(\Gamma_v)$ given by the fact above;
\item for $e\in E(S)$ take either $\lambda_e=\bbR^{n-1}$ or $\lambda_e=\Hat l_e$ (which coincides
with the axis of $\Gamma_e$ in $Y_{t(e)}$).
\end{itemize*}

The end completion being canonical, the group $\Gamma$ acts naturally on the disjoint union 
$\dunion_{v\in \cup V_F(S)} Y_v$, but we still have to define an action on
$\dunion_{v\in \Vab(S)} Y_v$ (the set of copies of $\bbR^{n-1}$). 

Remember that for $v\in V_F(S)$, the branching locus of axes in $Y_v$ is the set 
$$B_v=\bigcup_{\text{$e\neq e'\in E(S)$ incident on $v$}} \big[\lambda_e\cap \lambda_{e'}\big].$$
Since two distinct edges incident on $v$ have non-commuting fixators,
Fact \ref{fact_small} implies that $B_v$ is a countable union of sets of magnitude at most $n-2$.
Let $D_v=\{d(x,x')\,|\ x,x'\in B_v\}\subset\bbR^n$ be the set of mutual distances 
between points of the branching locus, and let $D=\bigcup \{D_v|v\in V_F(S)\}$.
As above, $D$ is a countable union of sets of magnitude at most $n-2$ of $\bbR^n$ (as a $\bbR^n$--metric space).

To get an action of $\Gamma$ on $\dunion_{v\in \Vab(S)} Y_v$, we just need to define 
an action of $\Gamma_v$ on $Y_v$ for one vertex $v$ of each orbit in $\Vab(S)/\Gamma$, and 
to extend this action equivariantly.
Now remember that for $v\in \Vab(S)$, $\Gamma_v=N_v\oplus \olG_v$, and $N_v$ comes with a natural
action on its axis in $T$ (and its translation length is infinitesimal), so we take an action of $N_v$ on $Y_v\simeq\bbR^{n-1}$
having the same translation length as in $T$. 
We choose the action of $\olG_v$ on $Y_v$ 
so that the translation length of any element of $\Gamma_v\setminus N_v$
lies in $\bbR^{n-1}\setminus D$.

We define the gluing maps $\phi_e$ inductively in a generic set as follows. 
First, up to changing $e$ to $\ol e$, we can assume that $t(e)\in \Vab(S)$ while $o(e)\in V_F(S)$.
For the first orbit of edges, we choose any $\Gamma_e$--equivariant gluing isometry $\phi_e\co \lambda_{\ol e}\ra \lambda_e$, 
and we extend this choice equivariantly.
Then, if some choices of gluing maps are already made
for some other edges incident on $v$, 
we choose $\phi_e$ so that
$\phi_e(B_{o(e)})$ does not meet $\phi_{e'}(B_{o(e')})$ for any edge $e'$ such that $t(e')=t(e)$ and on which $\phi_{e'}$
was already defined. This is possible since one can compose $\phi_e$ by any translation in $\bbR^{n-1}$,
and there are only countably many classes of translations mod $\bbR^{n-2}$ which are prohibited.
Then, we extend this choice equivariantly on the orbit of $e$.
This completes the definition of our graph of actions on $\bbR^{n}$--trees $\calg=(S,(Y_v)_{v\in V(S)},(\phi_e)_{e\in E(S)})$.

To sum up, our generic choices with respect to the branching locus of axes ensure that the following holds:
\begin{lem}\label{resume_phi}
\begin{itemize*}
\item For any vertex $v\in \Vab(S)$, the translation length of any element of $\Gamma_v\setminus N_v$
in $Y_v$ does not lie in $D$.
\item Assume that two edges $e,e'$ are incident on a common vertex $v\in \Vab(S)$ so
that $\phi_e\m B_{o(e)}\cap\phi_{e'}\m B_{o(e')}\neq \es$; then $e$ and $e'$ are in the same orbit.
\end{itemize*}
\end{lem}

\subsubsection{The dual action is free}

Let $T'=T_\calg$ be the $\bbR^n$--tree dual to the graph of actions defined above.
By lemma \ref{lem_flawless2}, the action of $\Gamma$ on $T'$ is flawless.
Therefore, Theorem \ref{thm_cyclic_flawless} will be proved
as soon as we prove that the action of $\Gamma$ on $T'$ is free.

\begin{lem}
The action of $\Gamma$ on $T'=T_\calg$ is free.
\end{lem}

\begin{proof}
We prove the freeness using our criterion (lemma \ref{lem_criterion_free}).
We will prove that no point of $Y_v$ gets identified with a point of $Y_w$ where
$v$ and $w$ are at distance more than $4$ in $S$.

Let's first consider a point $x\in Y_v$ where $v\in V_F(S)$ such that
$x$ is identified with two points $y_1\in Y_{u_1}$, $y_2\in Y_{u_2}$ for two distinct neighbours $u_1$, $u_2$ of $v$.
Then there are two edges $e\neq e'$  incident on $v$ such that $x\in\lambda_e\cap\lambda_{e'}$,
so $x$ lies in the branching locus of axes $B_v$.

Assume that there are two points of $Y_u$ and $Y_v$ which are identified in $T_\calg$,
for some $u,v$ at distance $5$ in $S$. 
Since every edge of $S$ joins a vertex in $\Vab(S)$ to a vertex in
$V_F(S)$,
there is a sub-path $(v_1,w,v_2)$ of $[u,v]\setminus\{u,v\}$ and some points
$x_1\in Y_{v_1}$, $y\in Y_w$, $x_2\in Y_{v_2}$
with $v_1,v_2\in V_F(S)$, $w\in \Vab(S)$. In particular,
for both $i=1,2$, $x_i\in B_{v_i}$.
By the second item of lemma \ref{resume_phi}, the edges $e_1=[v_1,w]$ and $e_2=[v_2,w]$
are in the same orbit. So let $g\in\Gamma$ sending $e_1$ on $e_2$. Note that $g\in\Gamma_w$ since
$w$ and $v_i$ are not in the same orbit. Since $N_w$ fixes all the edges incident on $w$, 
$g\in\Gamma_w\setminus N_w$. Now, since $x_1\in B_{v_1}$, $g.x_1\in B_{v_2}$, and $d(g.x_1,x_2)\in D_{v_2}$.
Let $y=\phi_{e_1}(x_1)$, so $g.y=\phi_{e_2}(g.x_1)$ so $d(y,g.y)=d(x_2,g.x_1)\in D$, a contradiction
with the first item of lemma \ref{resume_phi}.
\end{proof}

\section{D\'evissage theorem and corollaries}\label{sec_devissage}

\begin{SauveCompteurs}{devissage_simple}\devissageSimple
\end{SauveCompteurs}

This is a consequence of the following detailed version.

\begin{thm}[D\'evissage theorem, detailed version]\label{cyclic_devissage}
If a finitely generated freely indecomposable group $\Gamma$ acts freely on an $\bbR^n$--tree, 
then $\Gamma$ can be written as the fundamental group of a finite graph of groups $\calg$
with cyclic edge groups, finitely generated vertex groups,
with 3 types of vertices named \emph{abelian}, \emph{surface} and \emph{infinitesimal}, and such that the following holds:
\begin{itemize*}
\item each edge is incident to exactly one infinitesimal vertex;
\item for each abelian vertex $v$, $\Gamma_v$ is abelian maximal in $\Gamma$, 
$\Gamma_v=\olG_v\oplus N_v$ where $\olG_v$ is a finitely generated (maybe trivial) free abelian group,
$N_v$ is maximal cyclic in $G$, and the image in $\Gamma_v$ of all incident edges coincide with $N_v$;
moreover, if $v\neq v'$ are distinct abelian vertices, then $\Gamma_v$ does not commute with any conjugate of $\Gamma_{v'}$;
\item for each surface vertex $v$, $\Gamma_v$ is the fundamental group of a surface $\S$ with boundary holding
an arational measured foliation; there is one edge for each boundary component of $\S$, and the image of its edge group in $\Gamma_v$
is conjugate to the fundamental group of the corresponding boundary component of $\S$; 
\item for each infinitesimal vertex $v$, $\Gamma_v$ acts freely on an $\bbR^{n-1}$--tree; moreover, any element $g\in\Gamma\setminus\{1\}$ commuting
with an element of $\Gamma_v\setminus\{1\}$ either belongs to $\Gamma_v$, or is conjugate into $\olG_w$
where $w$ is an abelian vertex neighbouring $v$.
\end{itemize*}
Finally, $\calg$ is $4$--acylindrical and any non-cyclic abelian subgroup of $\Gamma$ is conjugate into a vertex group.
\end{thm}

\begin{proof}
Using Theorem \ref{thm_cyclic_flawless}, consider a free action $\Gamma\actson T'$
such that the action on the $\bbR$--tree $\olT '$ obtained by killing infinitesimals has cyclic arc fixators.
The Theorem is then a direct consequence of the abelian d\'evissage (Proposition \ref{abelian_devissage}):
the fact that $\Gamma$ and the edge groups of $\calg$ are finitely generated implies that
vertex groups are finitely generated.
\end{proof}

\begin{rem*}
The d\'evissage theorem does not claim that the splitting is non-trivial.
This occurs if $\Gamma$ is abelian, if $\Gamma$ is the fundamental group of a surface
with empty boundary, or if $\Gamma$ acts freely on some $\bbR^{n-1}$--tree.

It  follows from the commutative transitivity of $\Gamma$
that any non-cyclic maximal abelian subgroup of an infinitesimal vertex group
is maximal abelian in $\Gamma$. 
 However, some edge groups may fail to be  maximal cyclic in $\Gamma$ for some edges incident on surface vertices.
\end{rem*}

The following corollary is due to Sela and Kharlampovich--Myasnikov for limit groups \cite{Sela_diophantine1,KhMy_irreducible1}.
\begin{SauveCompteurs}{cor_FP}%
\corFP  
\end{SauveCompteurs}

\begin{proof}
For $n=1$, all the statements of the corollary
follow from Rips theorem which claims that $\Gamma$ is a free product of
finitely generated abelian groups and fundamental groups of closed surfaces (see \cite{GLP1,BF_stable}).

For $n>1$, we argue by induction and assume that the corollary holds for smaller values of $n$.
The conclusion of the corollary is stable under free product
since any non-cyclic abelian subgroup of a free product is conjugate into a vertex group.
Thus one can assume that $\Gamma$ is freely indecomposable. 
Then
the d\'evissage theorem says that $\Gamma$ is the fundamental group of a finite graph of groups $\calg$
with cyclic edge groups, and vertex groups satisfy the corollary by induction hypothesis.
If this splitting of $\Gamma$ is trivial, then $\Gamma$ is a vertex group and we are done.

The finite presentation of vertex groups implies that $\Gamma$ is finitely presented.
Moreover, induction hypothesis shows that $\Gamma$
has a finite classifying space,
and the cohomological dimension of $X$ is clearly at most $\max(2,r)$
 (see for instance \cite[Proposition VIII.2.4 and Ex.8b in VIII.6]{Brown_cohomology}).

We have the following bound about Betti numbers:
$$b_1(\Gamma)\geq \sum_{v\in V(\calg)} b_1(\Gamma_v) +b_1(\calg) -\# E(\calg)$$
where $b_1(\calg)$ denotes the first Betti number of the graph underlying $\calg$.
Indeed, consider the graph of groups $\calg_0$ obtained from $\calg$
by replacing edge group by a trivial group so that $\pi_1\calg_0$ is a free product of the vertex groups
and of a free group of rank $b_1(\calg)$. Since edge groups of $\calg$ are cyclic,
 one obtains $\Gamma$ from $\pi_1\calg_0$ by adding one relation for each edge of $\calg$,
and the inequality follows.
Now since $b_1(\calg)-1=\#E(\calg)-\#V(\calg)$, one gets that
$b_1(\Gamma)-1\geq \sum_{v\in V(\calg)} (b_1(\Gamma_v)-1)$.

By induction hypothesis, each term in the sum is non-negative.
In particular, $b_1(\Gamma)\geq b_1(\Gamma_v)$ for all $v\in V(\calg)$.
Thus if some vertex group is non-cyclic then $b_1(\Gamma)\geq 2$;
but all vertex groups cannot be cyclic because of acylindricity.

The d\'evissage theorem claims that a non-cyclic abelian subgroup $A$ fixes a vertex $v$ in the 
Bass--Serre tree $S$ of $\calg$. 
Let's prove that there are finitely many conjugacy classes of non-cyclic maximal abelian subgroups.
Since edge stabilizers are cyclic, such a subgroup $A$ fixes exactly one point in $S$.
Since there are only finitely many orbits of vertices in $S$, there remains to prove that
for any vertex $v\in V(S)$, there are only finitely many $\Gamma$--conjugacy classes of abelian subgroups $A$ which fix $v$.
The induction hypothesis says that there are at most finitely many such subgroups up to conjugacy in $\Gamma_v$, 
and therefore in $\Gamma$.

\newcommand{\Ab}{\mathrm{Ab}}
Denote by $\Ab(\Gamma)$ the set of conjugacy classes of abelian subgroups of $\Gamma$.
The argument above shows that the 
natural map $\dunion_{v\in V(\calg)} \Ab(\Gamma_v)\ra \Ab(\Gamma)$ is onto.
Therefore, 
\begin{eqnarray*}
  \sum_{A\in \Ab(\Gamma)} (\Rk A-1)&\leq& \sum_{v\in V(\calg)}\sum_{A\in \Ab(\Gamma_{v})} (\Rk A-1)\\
&\leq& \sum_{v\in V(\calg)} (b_1(\Gamma_{v})-1)\\
&\leq& b_1(\Gamma)-1.
\end{eqnarray*}
This terminates the proof of the corollary.
\end{proof}

\begin{cor}\label{cor_principal}
Consider a freely indecomposable, non-abelian, finitely generated group having a free action on an $\bbR^n$--tree.
Then  $\Gamma$ has a non-trivial splitting 
which is \emph{principal} in the following sense:
either $\Gamma=A*_C B$ or $\Gamma=A*_C$  where $C$ is maximal abelian in $\Gamma$, 
or $\Gamma=A*_C (C\oplus \bbZ^k)$.
\end{cor}

\begin{proof}
We argue by induction on $n$. The statement is clear for $n=1$.
Otherwise, consider the graph of groups given by the d\'evissage theorem. If $\calg$ contains a surface-type vertex,
then cutting along an essential curve provides a splitting over a cyclic subgroup which is maximal abelian.
If $\calg$ contains an abelian-type vertex $v$, write $G_v=N_v\oplus \olG_v$. 
If $\olG_v$ is trivial, then $N_v=G_v$ is maximal abelian, so each of the edges of $\calg$ incident on $v$
provides a principal splitting of $\Gamma$.
If $\olG_v$ is non-trivial, then $\Gamma=A*_{N_v}(N_v\oplus \olG_v)$
where $A$ is the fundamental group of the graph of groups $\calg'$ obtained from $\calg$ by replacing
the vertex group $\Gamma_v=N_v\oplus \olG_v$ by the cyclic group $N_v$.
If $\calg$ has no abelian-type and no surface-type vertex, then $\calg$ consists in a single infinitesimal vertex.
This means that $\Gamma$ acts freely on an $\bbR^{n-1}$--tree.
\end{proof}


\end{document}